\documentclass[11pt,a4paper]{amsart}
\usepackage{amsmath,amssymb,amsthm,graphicx,mathrsfs,url,latexsym,enumerate}
\usepackage[utf8]{inputenc}
\usepackage{textcomp}
\usepackage[english]{babel}
\usepackage{dsfont}

\numberwithin{equation}{section}

\newtheorem{hyp}{Assumption}

\newtheorem{thm}{Theorem}[section]
\newtheorem{lem}[thm]{Lemma}
\newtheorem{slem}[thm]{Sub-lemma}
\newtheorem{prop}[thm]{Proposition}
\newtheorem{cor}[thm]{Corollary}
\newtheorem{rem}[thm]{Remark}

\newcommand{\be}{\begin{equation}}
\newcommand{\ee}{\end{equation}}
\newcommand{\ba}{\begin{aligned}}
\newcommand{\ea}{\end{aligned}}
\newcommand{\bp}{{\it Proof. }}
\newcommand{\Vol}{\operatorname{Vol}}
\newcommand{\diam}{\operatorname{diam}}
\newcommand{\ep}{\hfill $\square$\\}
\newcommand{\N}{\mathbb N}

\newcommand{\R}{\mathbb R}
\newcommand{\C}{\mathbb C}
\newcommand{\Z}{\mathbb Z}

\newcommand{\Id}{\operatorname{Id}}
\newcommand{\supp}{\operatorname{supp}}

\newcommand{\Span}{\operatorname{Span}}
\newcommand{\vh}{\mathbf{h}}

\def\<{\langle}
\def\>{\rangle}
\def\indic{\mathds{1}}

\def\bbb{{\mathcal B}}\def\ccc{{\mathcal C}}
\def\eee{{\mathcal E}}\def\fff{{\mathcal F}} 
 \def\lll{{\mathcal L}}
 \def\ooo{{\mathcal O}}

\def\vvv{{\mathcal V}}

\def\Br{{\mathscr B}}

\def\Vr{{\mathscr V}}

\title[Metropolis algorithm on singular domains]{Spectral asymptotics for Metropolis algorithm on singular domains}
\author{L. Michel}

\begin{document}
\begin{abstract} We study the Metropolis algorithm on a bounded connected domain $\Omega$ of the euclidean space with proposal kernel localized at a small scale $h>0$. We consider the case of a domain $\Omega$ that may have cusp singularities. For small values of the parameter $h$ we prove the existence of 
a spectral gap $g(h)$ and study the behavior of $g(h)$ when $h$ goes to zero. 
As a consequence, we obtain exponentially fast return to equilibrium in total variation distance.
\end{abstract}
\maketitle
 \tableofcontents
 
\section{Introduction}\label{sec1}

Let $\Omega $
be a bounded connected open subset of $\mathbb{R}^d$ and let $\rho(x)$ be a positive measurable function on $\bar{\Omega}$ such that 
\be\label{eq:bornerho}
\forall x\in\Omega,\;m\leq \rho(x)\leq M
\ee
for some constants $m,M>0$.
We denote 
$\mu_\rho=\rho(x)dx$ the associated measure on $\Omega$ and we assume that   $\mu_\rho(\Omega)=\int_\Omega \rho(x)dx=1$. 
We consider the Metropolis algorithm associated to the density $\rho$ defined as follows.
For all $h\in ]0,1]$, we define the distribution kernel
\be\label{eq:defin-convkernel}
k_{h,\rho} (x,y)=h^{-d} \phi(\frac{x-y}{h}) \min (\frac{\rho(y)}{\rho(x)},1)
\ee
where $\phi(z)=\frac{1}{\Vr_d} \mathds{1}_{B(0,1)}(z)$, $B(0,1)$ denotes  the open unit ball  in $\R^d$ and $\Vr_d$ is the volume of $B(0,1)$.
The Metropolis kernel, is then given by 
\be\label{eq:defin-metrokernel}
t_{h,\rho} (x,dy)=m_{h,\rho}(x)\delta_{y=x}+k_{h,\rho} (x,y)dy
\ee
where $m_{h,\rho}(x)=1-\int_\Omega k_{h,\rho}(x,y)dy$. The kernel $t_{h,\rho} (x,dy)$ is clearly a Markov kernel on the state space $\Omega$ and 
the associated operator 
\be\label{eq:defin-metroop}
T_{h,\rho}(u)(x)= m_{h,\rho}(x) u(x)+ \int_\Omega k_{h,\rho} (x,y) u(y)dy
\ee
is a Markov operator. Throughout the paper, we sometimes omit the dependance of this operator with respect to $\rho$ and write $T_h$ instead of 
$T_{h,\rho}$ when there is no ambiguity.
A straightforward  computation shows that  $T_{h,\rho}$ is self-adjoint on $L^2(\Omega, \rho(x)dx)$ which implies 
in particular that the measure $\mu_\rho$ is stationary for the kernel $t_{h,\rho}(x,dy)$. As a consequence, the iterated kernel 
$t_{h,\rho}^n(x,dy)$ converges to the measure $\mu_\rho$ as $n\rightarrow\infty$, which explains the use of this kernel to sample the measure 
$\mu_\rho$. 

Introduced in \cite{MRRTT} to compute thermodynamical functionals by Monte-Carlo method, this algorithm has shown an impressive efficiency
and is now used as a routine in many domains of science. From a theoretical point of view, the computation of the speed of convergence of the
algorithm aroused many works in the setting of discrete state spaces (see \cite{Di09}, \cite{DiSa98} for introduction to this topic and references).  In \cite{DiLeMi11}, we obtained first results on a continuous state space in the limit $h\rightarrow 0$. More precisely, given  a bounded domain $\Omega$ of $\R^d$ with Lipschitz boundary we proved that the operator $T_h$ admits a spectral gap $g(h)$ of order $ h^2$ and for smooth densities $\rho$, we did compute the limit of $h^{-2}g(h)$. Eventually, we obtained some total variation estimates 
\be\label{eq:star1}
\sup_{x\in\Omega}\Vert t_{h,\rho}^n(x,dy)-d\mu_\rho(y)\Vert_{TV}\leq Ce^{-ng(h)}
\ee
for some constant $C>0$ independent of $h$. In this approach the fact that 
$\partial \Omega$ has Lipschitz regularity plays a fundamental role at several stages. A natural question is then to explore situations where this regularity assumption on $\partial\Omega$ fails to be true. 
In the present paper, we consider the case where $\partial \Omega$ may have cuspidal singularities.
More precisely we introduce the following assumption:
\begin{hyp}\label{hyp1}
There exist a  finite collection of open subsets of $\R^d$, $(\omega_i)_{i\in I_r\cup I_c}$ such that 
$\partial \Omega\subset (\cup_{i\in I_r\cup I_c} \omega_i )$ and
\begin{itemize}
\item[i)] for all $i\in I_r$,  $\partial \Omega \cap \omega_i$ has  Lipschitz regularity, 
\item[ii)] for all $i\in I_c$, there exists a closed submanifold $S_i$ of $\R^d$ with dimension $d''_i$, and there exist $\alpha_i>1,r_i>0,\epsilon_i>0$ such that in the neighbourhood of every point $p\in S_i$ there exists a coordinate system  
$(x_1,x', x'')\in\mathbb{R}^d=\mathbb{R}\times \mathbb{R}^{d'_i}\times \mathbb{R}^{d''_i}$, in which $p=(0, 0, 0)$ and
\be\label{eq:formcusp}
\Omega\cap \omega_i = \{(x_1,x',x''),0<x_1<\epsilon_i,\, \lvert x'\lvert_{d'_i}<x_1^{\alpha_i},\, \lvert x''\lvert _{d''_i}<r_i\}
\ee

\noindent where $\lvert \cdot \lvert _k$ stands for the euclidean norm on $\mathbb{R}^k$.

\end{itemize}
\end{hyp}

\noindent
Throughout the paper we will denote 
\be\label{eq:defgamma}
\gamma=\max_{i \in I_c}(\alpha_i-1)d_i'.
\ee
In our main results we need the cusp singularities to be not too sharp. We then introduce the following
\begin{hyp}\label{hyp2}
The constant $\gamma$ defined by \eqref{eq:defgamma}
satisfies $0<\gamma<2$.
\end{hyp}
Observe that as soon as $I_c$ is non empty (that is there exists some cusps on the boundary), one has $\gamma>0$.
Under the above assumption one has the following rough localization of the spectrum of $T_h$. The proof of this result will be given in the next section.

\begin{prop}\label{prop:spec-ess-Th}
Assume that Assumption \ref{hyp1} holds true. Then there exist $\delta_1,\delta_2>0$ and $ h_0>0$ such that for all $h\in ]0,h_0], $ $\sigma(T_h)\subset [-1+\delta_1h^\gamma,1]$ and $\sigma_{ess}(T_h)\subset [-1+\delta_1h^\gamma,1-\delta_2h^\gamma]$
where $\gamma$ is defined by \eqref{eq:defgamma}.

\end{prop}

From the above result, it is clear that the spectrum of $T_h$ in the interval $[1-Ch^{\gamma},1]$ is made of eigenvalues of finite multiplicity. 
Our first main result will provide precise informations on the spectrum of $T_h$ in a box
 $[1-Ch^{2},1]$ under  smoothness assumptions on the density $\rho$. For $\rho\in \ccc^1(\overline\Omega)$, we introduce the associated diffusion operator $L_\rho$ defined in a weak sense as follows. Given $u\in H^1(\Omega)$, let $\ell_u:H^1(\Omega)\rightarrow\C$ be defined by
$$
\ell_u(v)=\int_\Omega{\nabla \bar u}\nabla v\,d\mu_\rho+\int_\Omega \bar u  v\,d\mu_\rho
$$
where we recall that $d\mu_\rho=\rho(x)dx$.
We define the domain of $L_\rho$ as the set of 
functions $u\in H^1$ such that $\ell_u$ is continuous for the $L^2$ topology:
$$
D(L_\rho )=\{u\in H^1(\Omega),\,\exists C_u>0,\;\forall v\in H^1(\Omega),\;\vert \ell_u(v)\vert\leq C_u\Vert v\Vert_{L^2}\}
$$
Observe that $D(L_\rho)$ is not empty since it contains $\ccc_c^\infty(\Omega)$ (here we use the fact that $\rho$ is $\ccc^1$).
Since $H^1(\Omega)$ is dense in $L^2(\Omega)$ then  for any $u\in D(L_\rho)$, $\ell_u$ can be extended as a continuous linear form on $L^2(\Omega)$ and by Riesz Theorem, there exists a unique $f\in L^2(\Omega)$ such that 
$$\ell_u(v)=\<f, v\>_{L^2(\rho)},\;\forall v\in H^1(\Omega).$$
We then set $L_\rho u=-u+f$. From Theorem 3.6 in \cite{He13}, we know that $D(L_\rho)$ is dense in $H^1(\Omega)$ and that 
$\Id+L_\rho:D(L_\rho)\rightarrow L^2(\Omega)$  is bijective with bounded inverse. 
Now, it follows from Assumption \ref{hyp1} and the Theorem of section 8.3 in \cite{MaPo97} that the embedding $H^1(\Omega)\hookrightarrow L^2(\Omega)$ is compact and hence the resolvent 
$(\Id+L_\rho)^{-1}$ is compact. We  introduce the sequence 
 $\nu_0<\nu_1 <\nu_2<...$ of the distinct eigenvalues of $L_\rho$ with associated multiplicities $m_j$. Since,  $L_\rho$ is clearly non-negative and 
 $0$ is a simple eigenvalue, it follows $\nu_0=0$ and $m_0=1$.

\begin{thm}\label{th:spec-reg} 
Suppose that $\rho\in \ccc^1(\bar\Omega)$ satisfies \eqref{eq:bornerho}. Suppose that Assumptions \ref{hyp1} and \ref{hyp2} are verified.
Let $R>0, \epsilon>0$ and $J>0$ such that for all $j\leq J$, $\nu_j<R$ and for all 
$j<J$, $\nu_{j+1}-\nu_j >2\epsilon$. Then there exists $h_0>0$ such that for all $h\in ]0,h_0]$,

\be\label{1.5}
\sigma({1-T_{h}\over h^2})\cap ]0,R] \subset \cup_{j\geq 1} [\nu_j-\varepsilon, \nu_j+\varepsilon]
,\ee
 and the number of eigenvalues of ${1-T_h\over h^2}$ counted with multiplicities, in the interval 
 $[\nu_j-\varepsilon, \nu_j+\varepsilon]$,  is equal to $m_j$. 
\end{thm}

Observe that this theorem is the analogous of Theorem 1.2 in \cite{DiLeMi11}. Here we assume Assumption \ref{hyp2} to insure that there is no essential spectrum in the interval $[1-Ch^2,1]$. The case where $\gamma\geq 2$ seems more difficult to deal with since in this case the
 eigenvalues would be embedded in the essential spectrum.

If we drop the smoothness assumption on the density $\rho$ we get the following results.

\begin{thm}\label{th:spec-rough} Assume that $\rho$ is a measurable function satisfying \eqref{eq:bornerho}. Suppose that Assumptions \ref{hyp1} and \ref{hyp2} are  verified. Let $\delta_1,\delta_2>0$ be as in Prop. \ref{prop:spec-ess-Th}.
There exists $C,h_0>0$ 
such that for  any $h\in ]0,h_0]$,
 the following hold true:

\begin{enumerate}
\item[i)] The spectrum $\sigma(T_h)$ of $T_{h}$ is contained in  $[-1+\delta_1h^\gamma,1]$, 
$1$ is a simple eigenvalue of $T_{h}$, 
and $\sigma(T_h)\cap [1-\delta_2h^\gamma,1]$ is discrete. 
\item[ii)] The spectral gap $g(h):={\operatorname{dist}}(1,\sigma(T_h)\setminus\{1\})$ satisfies 
\be\label{gap3}
\frac 1Ch^2 \leq g(h) \leq C h^2.
\ee
\end{enumerate}
\end{thm}
As we shall see later, using \eqref{eq:bornerho} and comparaison of Dirichlet forms, this theorem is essentially  a consequence of Theorem \ref{th:spec-reg}.
From this spectral result we deduce estimates on the speed of convergence of the iterated kernel $( t^n_{h,\rho}(x, dy))$ towards the stationary measure $\mu_\rho$. We recall that the total variation distance between two probability measures $\mu$ and $\nu$ is defined by 
$\|\mu-\nu\|_{TV}=\sup_{A\in\Br}|\mu(A)-\nu(A)|$ where $\Br$ denotes the set of Borel set. Moreover, one has the following
\be\label{eq:TVdist1}
\|\mu-\nu\|_{TV}=
{1\over 2}\sup_{f\in L^\infty, \vert f\vert\leq 1}\vert \int fd\mu-\int fd\nu\vert.
\ee
\begin{thm}\label{th:TVestim}
Assume that $\rho$ is a measurable function satisfying \eqref{eq:bornerho}. Suppose that Assumptions \ref{hyp1} and \ref{hyp2} are is verified.
There exists $C,h_0>0$ 
such that for  any $h\in ]0,h_0]$ one has 
\be\label{1.7}
sup_{x\in \Omega}\Vert t^n_{h,\rho}(x, dy)-\mu_\rho\Vert_{TV} \leq C h^{-\gamma-\frac d2}e^{-n g(h)(1+O(h^{2-\gamma})) }.
\ee
for all $n\in\N$.
\end{thm}
Compare to Theorem 1.1 in \cite {DiLeMi11}, the estimate \eqref{eq:TVdist1} above suffers a loss of $h^{-\frac d 2-\gamma}$ in front of the exponential.
This loss is the natural loss when you go from convergence in $L^2$ sense (which follows from the spectral gap) to convergence in total variation. In 
\cite{DiLeMi11}, we used sophisticated tools (Nash estimates, Weyl asymptotics) to absorb this loss. In the present case, this strategy fails because of the cusp where nice estimates of eigenfunctions of $T_h$ can not be obtained from (see Lemma \ref{lem:decompHF}). However, let us emphasize that this prefactor implies only a logarithmic loss in the time needed to reach equilibrium ($h^{-2}\log(h)$ instead of 
$h^{-2}$).

The proof of the above theorems follows the general strategy of \cite{DiLeMi11}. 
In section \ref{sec1}, we prove Proposition \ref{prop:spec-ess-Th}. In order to prove Theorem \ref{th:spec-reg}   one uses minimax principle
and quasimodes built from the eigenfunctions of $L_\rho$  to prove that $h^{-2}(1-T_h)$ has at least $m_j$ eigenvalue near $\nu_j$. The proof of the converse inequality is more difficult and requires to prove some regularity property of eigenfunctions of $1-T_h$. This is done by mean of a dyadic decomposition of the cusp in section \ref{sec3}. Using  these constructions we prove the main theorems in section \ref{sec4}. In a separate appendix we prove a gluing lemma of $H^1$ functions which is crucially used in the proof of the main result.

We conclude this introduction with some notations used in the sequel.
On $\R^d$, we will denote by $|x|_d$ the euclidean norm of a vector $x$. When there is no ambiguity we will drop the index $d$ and simply write 
$|x|$.
Given a function $f:x=(x_1,x',x'')\in\R^{1+d'+d''}\rightarrow f(x_1,x',x'')\in\R$ we will denote by $\nabla'f(x)\in \R^{d'}$ (resp. $\nabla''f(x)\in \R^{d''}$)  the gradient of $f$ in the $x'$ variable (resp. $x''$ variable).
Given two quantities $u_t,v_t$ depending on a parameter $t$, we denote $u\asymp v$ if  there exists $C>0$ such that $\frac1 C u_t \leq v_t\leq Cv_t$
for all $t$.\\

\emph{Acknowledgment:} The author would like to thank  warmly G. Lebeau for numerous fruitful discussions  on this work.
The author is member of the ANR project QuAMProcs 19-CE40-0010-01.

\section{Rough localization of the spectrum }\label{sec2}
In this section, we give a proof of Proposition \ref{prop:spec-ess-Th}.
We first show that the operator $K_{h,\rho}: f\mapsto \int_\Omega k_{h,\rho}(x,y) f(y)\rho(y)dy$ is compact on $L^2(\Omega,\rho(y)dy)$.
Let $(\phi_n)$ be a sequence of continuous functions such that $\phi\leq\phi_n\leq 1$ and $( \phi_n)$ converges to $\phi$ in $L^2(\rho(x)dx)$ when $n\rightarrow \infty$. Consider the sequence of kernels $k_{n,h,\rho}=h^{-d} \phi_n(\frac{x-y}{h}) \min (\frac{\rho(y)}{\rho(x)},1)$ and let $K_{n,h,\rho}$ be the associated operators. Then 
$(K_{n,h,\rho})$ converges to $K_{h,\rho}$ in $\lll(L^2, L^2)$ when $n\rightarrow\infty$. 
%
On the other hand, since the kernels $k_{n,h,\rho}$ are continuous, the operators $(K_{n,h,\rho})$ are compact
and hence $K_{h,\rho}$ is compact.

Let us prove that $\sigma_{ess}(T_h)\subset [-1,1-Ch^\gamma]$. Thanks to Weyl criterium and compactness of $K_{h,\rho}$
 it is sufficient to prove that $\sup_{x\in \Omega} m_{h,\rho}(x)\leq 1-Ch^{\gamma}$.
Since
$$1-m_{h,\rho}(x)\geq \frac{mh^{-d}}{M \Vr_d}\int_\Omega \mathds{1}_{\lvert x-y\lvert <h} dy,$$
with $m,M$ given by \eqref{eq:bornerho}
the proof reduces to show that there exists $C,h_0>0$ such that 
\be\label{eq:minortheta}
\forall h\in]0,h_0],\;\forall x\in\Omega,\;\theta_h(x)\geq Ch^{d+\gamma}
\ee
where 
$\theta_h(x):=\int_\Omega \mathds{1}_{\lvert x-y\lvert <h} dy$.
Consider the family of subsets $\omega_i$ of 
Assumption \ref{hyp1} and let $\ooo_i=\Omega\cap\omega_i$. By a compactness argument, we can assume that there exists a family of open sets $(\tilde \omega'_i)$, such that 
$ \bar \omega'_i\subset \omega_i$ for all 
$i\in I_r\cup I_c$ and Assumption \ref{hyp1} holds true with the $\omega'_i$. It follows that 
$\Omega=\cup_{i\in J}^m\ooo'_i$ with $\ooo'_i= \omega'_i\cap\Omega$ 
where $J=I_r\cup I_c\cup\{0\}$, and $\omega'_0$ is an open subset of $\Omega$ such that $d(\bar\omega'_0,\partial \Omega)>0$.
Let us now estimate the function $\theta_h$ on each $\ooo'_i$.

%

We first observe that for $0<h_0<d(\bar\omega'_0,\partial \Omega)$ and $h\in]0,h_0]$, one has $B(x,h)\subset\Omega$ for any $x\in \ooo'_0$ and hence $\theta_h(x)=h^d\Vr_d$ which establishes the bound \eqref{eq:minortheta} on $\ooo'_0$.
Let us now study $\theta_h$ on $\ooo'_i$, $i\in I_r\cup I_c$. Taking $h_0>0$ sufficiently small, we can assume that for
all $h\in]0,h_0]$ one has $\omega'_i+B(0,h)\subset \omega_i$ for all $i\in I_r\cup I_c$. Hence,
if $\varphi:U_{i}\rightarrow \omega_{i}$ is a smooth local change of coordinates then
for any $x\in \ooo'_i$, one has
\be\label{eq:minormh1}
\begin{split}
\theta_h(x)&=\int_\Omega \mathds{1}_{\lvert x-y\lvert <h} dy\geq \int_{\ooo_{i}} \mathds{1}_{\lvert \varphi(\varphi^{-1}(x))-y\lvert <h} dy\\
&=\int_{U_{i}^+}J_\varphi(y) \mathds{1}_{\lvert \varphi(\varphi^{-1}(x))-\varphi(y)\lvert <h} dy
\end{split}
\ee
where  $J_\varphi(y)$ denotes the Jacobian of $\varphi$ and $U_{i}^+=\varphi^{-1}(\ooo_i)$. On the other hand, since $\varphi$ is a smooth function, there exists $C>0$ such that for all $ u,v\in U_i$, 
$\vert\varphi(u)-\varphi(v)\vert\leq C\vert u-v\vert$.
Combined with \eqref{eq:minormh1}, this implies
\be\label{eq:minormh2}
\theta_h(x)\geq \int_{U_{i}^+}J_\varphi(y) \mathds{1}_{\lvert \varphi^{-1}(x)-y\lvert <h/C} dy\geq \tilde C\int_{U_{i}^+}\mathds{1}_{\lvert \varphi^{-1}(x)-y\lvert <h/C}dy
\ee
for some positive constant $\tilde C$ such that $| J_\varphi|\geq \tilde C$ on $U_{i}$. This minoration shows that in order to get some lower bound on $\theta_h$, we can suppose that we are in any suitable system of coordinates.

Suppose that $i\in I_r$. By a  Lipschitz change of coordinates it is shown in \cite{DiLeMi11} that there exists some  constants $c_1,c_2>0$ such that 
\be\label{eq:minormh2}
\theta_h(x)\geq c_1\int_{x_1\geq 0}\mathds{1}_{\lvert x-y\lvert <h} dy\geq c_2 h^d
\ee
 for all $x\in\ooo'_i$. Combined with the definition of $\theta_h$, this shows that 
$(1-m_h(x))\geq c_3$ for some $c_3>0$ independent of $h$.

Suppose now that $i\in I_c$  and that  $\omega_{i}$ is like in ii) of Assumption \ref{hyp1}. Using  a suitable change of coordinates, we can assume
that  
there exist $\alpha>1, r>0,\epsilon>0$ such that
$$\ooo_i=\Omega\cap \omega_{i}= \{(x_1,x',x''),0<x_1<\epsilon,\, \lvert x'\lvert_{d'}<x_1^{\alpha},\, \lvert x''\lvert _{d''}<r\},$$
where $d'$, $d''$ are the local dimension appearing in Assumption \ref{hyp1} whose dependance with respect to the index $i$ is omitted.
Moreover, we can also assume that $\ooo'_i=\ooo_i\cap\{0<x_1<\epsilon/2\}\cap\{|x''|<r/2\}$.
Endowing  $\ooo_{i}$ with the equivalent norm
$$ \lvert (x_1,x',x'')\lvert_\infty =\max \{ \lvert x_1\lvert, \lvert x'\lvert _{d'}, \lvert x''\lvert_{d''} \},$$
 it is sufficient to find a lower bound for $\int_\Omega \mathds{1}_{ \lvert x-y \lvert _\infty <h} dx$  when $x$ varies in  $\ooo'_i$. For such $x$, one has
\begin{equation}\label{eq:minorcusp}
\begin{split}
\int_\Omega \mathds{1}_{ \lvert x-y\lvert_\infty <h} dy
&=\int_\Omega \mathds{1}_{ \lvert (x_1,x',x'')-(y_1,y',y'')\lvert _\infty<h}dy_1dy'dy''\\
&=\int_{\lvert y''\lvert_{d''} <r, \lvert y'\lvert_{d'}<y_1^\alpha, 0<y_1<\epsilon} \mathds{1}_{\lvert x''-y''\lvert _{d''}<h}\\
&\phantom{************} \mathds{1}_{\lvert x'-y'\lvert_{d'}<h} \mathds{1}_{\lvert x_1-y_1\lvert <h}dy_1dy'dy''\\
& \geq ch^{d''} W_h(x,x')
\end{split}\end{equation}
where $c$ is a positive constant and 
\be\label{eq:defWh}
W_h(x,x'):= \int_{\lvert y'\lvert_{d'}<y_1^\alpha, 0<y_1<\epsilon}\mathds{1}_{\lvert y'-x'\lvert_{d'}<h}\mathds{1}_{\lvert x_1-y_1\lvert <h} dy'dy_1.
\ee 
Denoting
$$\ccc=\{(y_1,y')\in \R\times\R^{d'},\,\vert y'\vert_{d'}<y_1^\alpha, \;0<y_1<\epsilon\}$$ and 
$$
D_h(x_1,x')=\{y\in \ccc,\; \lvert y'-x'\lvert_{d'}<h\text{ and }\lvert x_1-y_1\lvert <h\},
$$
we have $W_h(x_1,x')=vol(D_h(x_1,x'))$ and thanks to \eqref{eq:minorcusp}, one has to prove that 
$W_h(x_1,x')\geq ch^{\alpha d'+1}$ for some uniform constant $c>0$.
We first observe that it holds true for $(x_1,x')=(0,0)$, since one has (using $\alpha>1$)
\begin{equation}\label{eq:minorcusp1}
\begin{split}
W_h(0,0)&= \int_0^h\int_{\R^{d'}} \mathds{1}_{B_{d'}(0,y_1^\alpha)}(y')dy'dy_1
=\Vr_{d'}\int_0^hy_1^{\alpha d'}dy_1=c h^{\alpha d'+1}.
\end{split}\end{equation}
We now decompose the cusp into three zones  that we treat differently: $\{0<x_1\leq h/2\}$, $\{h/2<x_1<(\delta h)^{\frac 1 \alpha}\}$ and
$\{(\delta h)^{\frac 1 \alpha}<x_1<\epsilon\}$, where $\delta>0$ will be chosen sufficiently small. 
\begin{itemize}
\item Suppose first that  $x_1 < \frac h 2 $, then since $\alpha>1$, one has for $h$ small enough
$
D_h(x_1,x')\supset\{(y_1,y'),\;\vert y_1\vert<h/2\}\cap \ccc.
$
Combined with \eqref{eq:minorcusp1}, this yields
$$
W_h(x_1,x')\geq\int_{\lvert y'\lvert_{d'}<y_1^\alpha}
\mathds{1}_{\lvert y_1\lvert <h/2} dy_1dy'=W_{\frac h 2}(0,0)=c h^{\alpha d'+1}
$$
which is the required lower bound on $W_h$.
\item Suppose now that  $h/2\leq x_1 <(\delta h)^{\frac 1\alpha} $, then  
$$
D_h(x_1,x')\supset\{x_1<y_1<x_1+h,\;\vert y'\vert<x_1^\alpha\}.
$$
Indeed, if $\vert y'\vert<x_1^\alpha$ and  $y_1>x_1$ one gets immediately  $(y_1,y')\in \ccc$ and since $|x'|<x_1^\alpha<\delta h$ then 
$|x'-y'|\leq |x'|+|y'|\leq 2\delta h<h$ for $0<\delta<\frac 12$.
From the above inclusion, it follows 
$$
W_h(x_1,x')\geq  \Vr_{d'}x_1^{\alpha d'}h\geq  2^{-\alpha d'}\Vr_{d'}h^{\alpha d'+1}.
$$
\item Eventually, suppose that $(\delta h)^{\frac 1\alpha} \leq x_1<\epsilon$. We observe that the application $x'\mapsto W_h(x_1,x')$ is radial.
Hence, it suffices to estimate from below the application $t\in[0,x_1^\alpha]\mapsto W_h(x_1,x'_t)$ with $x'_t=(t,0,\ldots,0)$.

\begin{itemize}
\item[-] If $|t|<\delta h/2$ then  the inclusion 
$$
D_h(x_1,x_t')\supset\{x_1<y_1<x_1+h,\;\vert x_t'-y'\vert<\delta h/2\}
$$
implies 
$W_h(x_1,x_t')\geq h(\delta h/2)^{d'}=ch^{d'+1} $. 

\item[-]
If $\delta h/2\leq |t|<x_1^\alpha$ then
$$\{|y'-x'_{t_h}|<\delta h/4\}\subset \{|y'|<|t|\}\subset \{|y'|<x_1^\alpha\}$$
where $t_h=t-\delta h/4$ and $x'_{t_h}=(t_h,0,\ldots,0)$. Hence 
$$
D_h(x_1,x_t')\supset\{x_1<y_1<x_1+h,\;\vert x_{t_h}'-y'\vert<\delta h/4\}
$$
which implies again that $W_h(x_1,x_t')\geq ch^{d'+1} $. 
\end{itemize}
Summing up the above discussion, we have proved that  for any $i\in I_c$, $W_h(x_1,x')\geq ch^{\alpha_i d'+1}$ uniformly on 
$\ooo_i'$. Combined with \eqref{eq:minorcusp}, this proves that $\theta_h(x)\geq ch^{(\alpha_i-1)d'}$ uniformly on $\ooo_i'$.
\end{itemize}
Since the boundary of $\Omega$ is compact,  it follows from the above computations that there exists $c>0$ such that for all $x\in \Omega$ and all $h<h_0$,
$m_h(x)\leq 1-Ch^\gamma$
with $\gamma$ given by \eqref{eq:defgamma}. 
This proves that $\sigma_{ess}(T_h)\subset [-1,1-Ch^\gamma]$.

We now prove that $\sigma(T_h)\subset [-1+Ch^\gamma,1]$ which is equivalent  to show that 
$$\langle u+T_hu,u\rangle_{L^2(\rho)} \geq Ch^\gamma\lvert \lvert u\lvert \lvert ^2_{L^2(\rho)}$$
for all $u\in L^2(\Omega)$. For this purpose, we observe that thanks to \cite{DiLeMi11}, eq. (2.7), one has
$$\langle u+T_hu,u\rangle_{L^2(\rho)}\geq \frac{1}{2} \int_{\Omega\times \Omega} k_{h,\rho}(x,y)\lvert u(x)+u(y)\lvert ^2 \rho(x) dxdy.$$
Hence, it is sufficient to prove that there exist $C_0,h_0>0$ such that the following inequality holds true for all $h\in ]0,h_0]$ and all $u\in L^2(\Omega)$:

$$
 \int_{\Omega\times \Omega} k_{h,\rho}(x,y)\lvert u(x)+u(y)\lvert ^2 \rho(x) dxdy\geq C h^{\gamma} \Vert u\Vert^2_{L^2(\rho)}.
$$
Since $\rho$ is bounded from below, we can assume without loss of generality that $\rho=1$.
Following \cite{DiLeMi11}, we introduce a covering $(\nu_j)_j$ of $\Omega$ with $\nu_j\subset \Omega$ such that $\diam(\nu_j)<h$ and for some $C_1>0$ independent of $h$, the number of indices $k$ such that $\nu_j \cap \nu_k \neq \emptyset$ is less than $C_1$. Moreover,
since $\inf_\Omega \theta_h\geq Ch^{d+\gamma}$, we can also assume that there exists  a constant $C_2>0$ such that $\Vol(\nu_j)\geq C_2 h^{d+\gamma}$  for any $j$.
Then, we get as in \cite{DiLeMi11}
\begin{equation*}\begin{split}
C_1 \int_{\Omega\times \Omega} h^{-d} \phi(\frac{x-y}{h})&\lvert u(x)+u(y)\lvert^2dxdy\\
&\geq \sum_j \int_{\nu_j\times \nu_j} h^{-d}  \phi(\frac{x-y}{h}) \lvert u(x)+u(y)\lvert ^2 dxdy\\
& \geq \sum_j h^{-d} \frac{1}{\Vr_d }\int_{\nu_j\times \nu_j} \lvert u(x)+u(y)\lvert ^2 dxdy\\
&\geq \sum_j 2h^{-d} \frac{1}{\Vr_d } \Vol(\nu_j) \lvert \lvert u\lvert \lvert^2_{L^2(\nu_j)}\\
&\geq \frac{2C_2 h^\gamma}{\Vr_d } \lvert \lvert u\lvert \lvert ^2_{L^2(\Omega)}.
\end{split}\end{equation*}
This implies, 
$\langle u+T_hu,u\rangle \geq \tilde{C} h^\gamma \lvert \lvert u\lvert \lvert^2$
and finally
$\sigma(T_h)\subset [-1+Ch^\gamma, 1]. $ The proof of Proposition \ref{prop:spec-ess-Th} is complete.

\section{Regularity of eigenfunctions}\label{sec3}

The aim of this section is to prove regularity properties on families of eigenfunctions of $T_h$ associated to eigenvalues in 
$[1-Ch^\gamma,1]$.
Let us introduce the Dirichlet form of $T_h$
\be\label{eq:defDirich}
\bbb_{h,\rho}(f,g):=\<(1-T_h)f,g\>_{L^2(\rho)}
\ee
and $\eee_{h,\rho}(f)=\bbb_{h,\rho}(f,f)$.
One has
$$
\bbb_{h,\rho}(f,g)=\frac 1{2}\int_{\Omega\times \Omega} k_{h,\rho}(x,y)( f(x)-f(y))\overline{(g(x)-g(y))} \rho(x) dxdy
$$
and denoting $d\mu_\rho^2=\min(\rho(x),\rho(y)) dxdy$ we get 
$$
\bbb_{h,\rho}(f,g)=\frac 1{2 h^{d}\Vr_d}\int_{\Omega\times \Omega}\indic_{\vert x-y\vert<h}( f(x)-f(y))\overline{(g(x)-g(y))} d\mu_\rho^2(x,y).
$$
In particular, one has 
\be\label{eq:defDirich1}
\eee_{h,\rho}(f)=\frac 1{2 h^{d}\Vr_d}\int_{\Omega\times \Omega}\indic_{\vert x-y\vert<h}|f(x)-f(y)|^2 d\mu_\rho^2(x,y).
\ee
As mentioned before, we will sometimes drop index $\rho$ in the notations when it is unambiguous.
The following decomposition lemma is the key point in our analysis. 
\begin{lem}\label{lem:decompHF}
Let $(f_h)_{h\in]0,1]}$ be a family of function in $L^2(\Omega)$ such that $\Vert f_h\Vert _{L^2}\leq 1$ and
$\eee_h(f_h)\leq h^2$.
Then, there exists $C,C_0,h_0>0$ such that for all $h\in ]0,h_0]$, one has a decomposition
$f_h=f_{h,\ccc}+f_{h,L}+f_{h,H}$
with
\begin{itemize}
\item[-] $\supp(f_{h,\ccc})\subset\Gamma_{2h}$ with $\Gamma_h=\cup_{i\in I_c}\{x\in\Omega,\;d(x,S_i)<Ch^{\frac 1\alpha_i}\}$, $\alpha_i, S_i$ given by 
Assumption \ref{hyp1}
\item[-] $f_{h,L}$ and $f_{h,H}$ are supported in $\Omega\setminus \Gamma_h$ and 
$$\lvert \lvert \nabla f_{h,L}\rvert \rvert _{L^2}\leq C_0\text{ and }\lvert \lvert f_{h,H}\rvert \rvert _{L^2}\leq C_0 h$$

\end{itemize}

\end{lem}

This lemma is inspired from Lemma 2.2 in  \cite{DiLeMi11}. However, due to  the presence of cusps there is an additional term in the decomposition of 
$f_h$ for which we do not have nice estimates. Moreover, we have to face important complications in the proof. The next section is devoted to the proof of this lemma in the particular case where $\Omega$ is a model cusp.

\subsection{A model case}

In this section we consider the   case where the domain $\Omega$ is an exact cusp
\be\label{eq:defOmegamodel}
\Omega=\{(x_1,x',x''),\;0<x_1<1, \lvert x'\lvert _{d'}<x_1^{\alpha}, \lvert x''\lvert_{d''} <1\}.
\ee
Since there is no ambiguity, $\Omega$ denotes the above domain in this section and a general domain in the rest of the paper.
Since $\rho$ is bounded from below and above by positive constant, we can assume that $\rho=1$ 
without modifying the assumption $\eee_h(f_h)= \ooo(h^2)$.
One defines a dyadic partition $(\Omega_k)_{k\geq 0}$ of ${\Omega}$ in the following way:
$$\Omega_k:={\Omega}\cap \{\frac{1}{2^{k+1}} < x_1 < \frac{1}{2^k}\}, k\in\N.$$
For every $k\geq 0$, we define a change of variables 
\be\label{eq:deftauk}
\begin{array}{cc}
\tau_k:&\Omega_k\,\,\,\rightarrow\,\,\, \Omega_0\phantom{********}\\
& (x_1,x',x'')\mapsto (2^k x_1, 2^{k\alpha}x',x'')
\end{array}
\ee
whose jacobian is $j_k:=\det d\tau_k=2^{k(\alpha d'+1)}$. We also introduce the change of variable
\be\label{eq:def-hattauk}
\begin{array}{cc}
\hat \tau_k:&\Omega_k\,\,\,\rightarrow\,\,\,\Omega_1\phantom{********}\\
& (x_1,x',x'')\mapsto (2^{k-1} x_1, 2^{(k-1)\alpha}x',x'')
\end{array}
\ee
 and we observe that $\tau_k=\tau_1\circ\hat\tau_k$.
\subsubsection{Sobolev space and dyadic decomposition of cusps}
Throughout the paper we will use the following notation. Given a set $B$ a function  $f\in H^1(B)$, and some parameters $h,\tilde h, \bar h>0$, we denote
\be\label{eq:notnormdyad}
N_{\bar h,\tilde h,h}(f,B)=\Big(\Vert \bar h\partial_1 f\Vert_{L^2(B)}^2+\Vert \tilde h \nabla'f\Vert^2_{L^2(B)}
+\Vert h\nabla'' f\Vert^2_{L^2(B)}\Big)^{\frac 12}.
\ee
In order to lighten the notation we introduce the parameter $\vh=(\bar h,\tilde h,h)$ and we will often write 
$$
N_{\vh}(f,B)=N_{\bar h,\tilde h,h}(f,B).
$$
The following lemma gives an expression of Sobolev norms for dyadic decomposition of the domain $\Omega$.
\begin{lem}\label{lem:decomp-dyad-Sob1}
Let $f\in L^2(\Omega)$, then 
$$
\Vert f\Vert_{L^2(\Omega)}^2=\sum_{k\in\N}2^{-k(\alpha d'+1)}\Vert f\circ\tau_k^{-1}\Vert_{L^2(\Omega_0)}^2
$$
If one assume additionally that $f\in H^1(\Omega)$, then 
\begin{equation*}
\begin{split}
\Vert \nabla f\Vert_{L^2( \Omega)}^2&=\sum_{k\in\N}2^{-k(\alpha d'+1)}\Big(\Vert 2^k\partial_1(f\circ\tau_k^{-1})\Vert_{L^2(\Omega_0)}^2\\
&\phantom{****}+\Vert 2^{k\alpha}\nabla'(f\circ\tau_k^{-1})\Vert_{L^2(\Omega_0)}^2
+\Vert \nabla''(f\circ\tau_k^{-1})\Vert_{L^2(\Omega_0)}^2\Big)\\
&=\sum_{k\in\N}2^{-k(\alpha d'+1)}N_{2^k,2^{k\alpha},1}(f\circ\tau_k^{-1}, \Omega_0)^2
\end{split}
\end{equation*}
\end{lem}
\bp
Use the  partition $\Omega=\cup_{k\in\N}\Omega_k$, the change of variable $\tau_k$ and the chain rule.
\ep
Let 
\be\label{eq:deftheta}
\begin{split}
\theta:\;\; \R_+^*&\times\R^{d'}\times\R^{d''}\rightarrow\R_+^*\times\R^{d'}\times\R^{d''}\\
&(x_1,x',x'')\mapsto (x_1, x_1^{-\alpha}x',x'')
\end{split}
\ee
 and consider the open sets
$B_j:=\{\frac 1{2^{j+1}}<x_1<\frac 1{2^j},\,\vert x'\vert<1,\vert x''\vert<1\}$. Observe that $\theta$ is a $C^1$ diffeomorphism from 
$\Omega_j$ onto $B_j$. 
Hence, the maps 
\be\label{eq:defsigmak}
\sigma_k=\theta\circ\tau_k:\Omega_k\rightarrow B_0
\ee 
and 
\be \label{eq:def-hatsigmak}
\hat\sigma_k=\theta\circ\hat\tau_k:\Omega_k\rightarrow B_1
\ee are also   $C^1$ diffeomorphisms. Moreover, one has 
$
\sigma_k=\check\sigma_1\circ\hat\sigma_k
$
where  
\be\label{eq:defchecksigma1}
\check\sigma_1=B_1\rightarrow B_0,\;
\check\sigma_1=\sigma_1\circ\theta^{-1}=\theta\circ\tau_1\circ\theta^{-1}.
\ee
The following lemma express $L^2$ and $H^1$ norm in terms of the dyadic decomposition.
\begin{lem}\label{lem:decomp-dyad-Sob2}
One has the following estimates 
\be\label{eq:decomp-dyad-Sob1} 
\Vert f\Vert_{L^2(\Omega)}^2\asymp\sum_{k\in\N}2^{-k(\alpha d'+1)}\Vert f\circ\sigma_k^{-1}\Vert_{L^2(B_0)}^2
\ee
for any $f\in L^2(\tilde\Omega)$ and 
\be\label{eq:decomp-dyad-Sob2}
\Vert \nabla f\Vert_{L^2( \Omega)}^2\asymp\sum_{k\in\N}2^{-k(\alpha d'+1)}N_{2^k,2^{k\alpha},1}(f\circ\sigma_k^{-1},B_0)^2
\ee
for any $f\in H^1(\Omega)$.
Conversely, assume that $(f_k)_{k\in\N}$ is a sequence of functions of $H^1(B_0)$  such that 
\be\label{eq:decomp-dyad-Sobfini}
\sum_{k\in\N}2^{-k(\alpha d'+1)}(\Vert f_k\Vert_{L^2(B_0)}^2+N_{2^k,2^{k\alpha},1}(f_k,B_0)^2)<\infty
\ee
and $(f_k)_{\vert x_1=\frac 12}=(f_{k+1})_{\vert x_1= 1}$, where $(f_k)_{\vert x_1=a}$ denotes the trace of the $H^1$ function $f_k$ on $\{x_1=a\}$. Then the function $f:=\sum_{k=0}^\infty\indic_{\Omega_k} f_k\circ\sigma_k$ belongs to 
$ H^1(\Omega)$. Moreover, for such functions, one has
$$\Vert f\Vert_{H^1(\Omega)}^2\asymp \sum_{k\in\N}2^{-k(\alpha d'+1)}(\Vert f_k\Vert_{L^2(B_0)}^2+N_{2^k,2^{k\alpha},1}(f_k,B_0)^2).$$
\end{lem}
\bp
For any $j\geq 0$, $\theta$ defines a change of variable from $\Omega_j$ onto $B_j$.
A standard computation shows that there exists $C>1$ such that 
$$
\frac 1 C\Vert f\circ\theta^{-1}\Vert_{L^2(B_0)}\leq \Vert f\Vert_{L^2(\Omega_0)}\leq C\Vert f\circ\theta^{-1}\Vert_{L^2(B_0)},
$$
and
\begin{equation*}
\begin{split}
\frac 1 C&\Vert \nabla'(f\circ\theta^{-1})\Vert_{L^2(B_0)}\leq \Vert \nabla' f\Vert_{L^2(\Omega_0)}\leq C\Vert \nabla'(f\circ\theta^{-1})\Vert_{L^2(B_0)}\\
\frac 1 C&\Vert \nabla''(f\circ\theta^{-1})\Vert_{L^2(B_0)}\leq \Vert \nabla'' f\Vert_{L^2(\Omega_0)}\leq C\Vert \nabla''(f\circ\theta^{-1})\Vert_{L^2(B_0)}\\
\frac 1 C&\Vert \partial_1(f\circ\theta^{-1})\Vert_{L^2(B_0)}\leq \Vert \partial_1f\Vert_{L^2(\Omega_0)}+
\Vert \nabla' f\Vert_{L^2(\Omega_0)}\\
\frac 1 C&\Vert \partial_1f\Vert_{L^2(\Omega_0)}\leq \Vert \partial_1(f\circ\theta^{-1})\Vert_{L^2(B_0)}+
\Vert \nabla'(f\circ\theta^{-1})\Vert_{L^2(B_0)}
\end{split}
\end{equation*}
Combining these estimates with Lemma \ref{lem:decomp-dyad-Sob1}, we obtain 
 \eqref{eq:decomp-dyad-Sob1} and \eqref{eq:decomp-dyad-Sob2}.
 Conversely, assume that $f\in L^2(\Omega)$ is such that 
\eqref{eq:decomp-dyad-Sobfini} holds true. In order to prove that $f\in H^1(\Omega)$, it suffices to show that $f$ has no jump 
at $x_1=2^{-k}$. This exactly the condition $(f_k)_{|x_1=\frac 12}=(f_{k+1})_{|x_1=1}$.
\ep

\begin{rem}\label{rem-norm-Sob}
If one splits the sums in the above lemma into even and odd terms, one gets
\begin{equation*}\label{eq:decomp-dyad-Sob-split1}
 \begin{split}
\Vert f\Vert_{L^2(\tilde\Omega)}^2\asymp\sum_{k\in\N}&4^{-k(\alpha d'+1)}\Vert f\circ\sigma_{2k}^{-1}\Vert_{L^2(B_0)}^2\\
&\phantom{*******}+\sum_{k\in\N}4^{-k(\alpha d'+1)}\Vert f\circ\sigma_{2k+1}^{-1}\Vert_{L^2(B_0)}^2
\end{split}
\end{equation*}
Using the identity $\sigma_{2k+1}=\sigma_1\circ\theta^{-1}\circ \hat\sigma_{2k+1}$ with $\hat\sigma_{2k+1}$ defined by \eqref{eq:def-hatsigmak}, \eqref{eq:def-hattauk} and the fact that 
$\theta\circ\sigma_1^{-1}$ is a diffeomorphism from $B_0$ onto $B_1$, we get 
\begin{equation*}\label{eq:decomp-dyad-Sob-split2} 
\begin{split}
\Vert f\Vert_{L^2(\tilde\Omega)}^2\asymp\sum_{k\in\N}4^{-k(\alpha d'+1)}\Vert &f\circ\sigma_{2k}^{-1}\Vert_{L^2(B_0)}^2\\
&+\sum_{k\in\N}4^{-k(\alpha d'+1)}\Vert f\circ\hat\sigma_{2k+1}^{-1}\Vert_{L^2(B_1)}^2
\end{split}
\end{equation*}
Similarly, we get the following identity for the norm of the gradient
\begin{equation*}\label{eq:decomp-dyad-Sob-split3}
\begin{split}
\Vert \nabla f\Vert_{L^2(\tilde \Omega)}^2&\asymp\sum_{k\in\N}4^{-k(\alpha d'+1)}N_{4^k,4^{k\alpha},1}(f\circ\sigma_{2k}^{-1},B_0)^2\\
&\phantom{********}
+ \sum_{k\in\N}4^{-k(\alpha d'+1)}N_{4^k,4^{k\alpha},1}(f\circ\hat \sigma_{2k+1}^{-1},B_1)^2.
\end{split}
\end{equation*}
\end{rem}

\subsubsection{Extension of the operator and comparison of Dirichlet forms}
Let $f \in L^2({\Omega})$ be such that $\Vert f\Vert_{L^2}=1$ and $\eee_h(f)=\ooo(h^2)$. 
We observe that since the norms $\vert.\vert$ and $\vert.\vert_\infty$ are equivalent, there exists a constant $C>0$ such that 
\be\label{eq:compar-dirich-norme}
\eee_{\infty,\frac h C}(f)\leq\eee_h(f)\leq\eee_{\infty, Ch}(f)
\ee
where
\be\label{eq:defEinfty}
\eee_{\infty,h }(f):=\frac 1{2 h^{d}\Vr_d}\int_{ \Omega\times \Omega}\indic_{\vert x-y\vert_{\infty}<h}\lvert f(x)-f(y)\lvert ^2 dxdy.
\ee
Thanks to \eqref{eq:compar-dirich-norme}, one has
\begin{equation}\label{eq:decomp-dirich-scales1}
\begin{split}
\eee_{C h}(f)&\geq \eee_{\infty,h}(f)=\frac{1}{2\Vr_d h^d}\int_{ \Omega\times  \Omega}\indic_{\lvert x-y\lvert _\infty <h}\lvert f(x)-f(y)\lvert^2 dx dy\\
&\geq \sum_{k\geq 0} \frac{1}{2\Vr_d h^d} \int_{\Omega_k\times  \Omega_k}\indic_{\lvert x-y\lvert _\infty <h} \lvert f(x)-f(y)\lvert ^2 dxdy\\
& = \sum_{k\geq 0} \frac{1}{2\Vr_d h^d} \int_{ \Omega_0\times \Omega_0 }
\indic_{ \lvert x_1-y_1\lvert <2^k h, \lvert x'-y'\lvert <2^{k\alpha}h,\lvert x''-y''\lvert <h}\\
&\phantom{********}\lvert f\circ \tau_k^{-1}(x)-f\circ \tau_k^{-1}(y)\lvert ^2 j_k^{-1}(x) j_k^{-1}(y) dxdy\\
& = \sum_{k\geq 0} \frac{2^{-k(1+\alpha d')}}{2\Vr_d h^{d''} \tilde{h}_k^{d'} \bar{h}_k} \int_{\Omega_0\times  \Omega_0}\indic_{\lvert x_1-y_1\lvert <\bar{h}_k, \lvert x'-y'\lvert < \tilde{h}_k,\lvert x''-y''\lvert <h}\\
&\phantom{************}
\lvert f\circ \tau_k^{-1}(x)-f\circ \tau_k^{-1}(y)\lvert ^2 dxdy
\end{split}
\end{equation}
 where $\bar{h}_k=2^kh, \tilde{h}_k= 2^{k\alpha}h$. 
 Given any domain $A\subset\R^d$, one then introduces the Dirichlet form defined on $L^2(A)$ by
$$
\eee_{\bar{h},\tilde{h},h}^A(g)=\frac{1}{2\Vr_d h^{d''}\tilde{h}^{d'}\bar{h}}
\int_{A\times A}\indic_{ \lvert x_1-y_1\lvert <\bar{h}, \lvert x'-y'\lvert <\tilde{h}, \lvert x''-y''\lvert <h}\lvert g(x)-g(y)\lvert ^2 dxdy.
$$
Then the last inequality in \eqref{eq:decomp-dirich-scales1} reads
\be\label{eq:decomp-dirich-scales2}
\eee_{C h}(f)\geq \sum_{k\geq 0}2^{-k(1+\alpha d')}\eee^{ \Omega_0}_{\bar{h}_k,\tilde{h}_k,h}(f\circ \tau_k^{-1}).
\ee
The next step in the computation is to compare the Dirichlet form $\eee^{\Omega_0}_{\bar{h}_k,\tilde h_k,h}$ and
$\eee^{B_0}_{\bar{h}_k,\tilde h_k,h}$ associated respectively to the domains
$ \Omega_0$ and $B_0$. As a preliminary step, we need the following result.
\begin{lem}\label{lem:mult-param-dirich} Let $A$ be any open subset of $\R^d$ with Lipschitz boundary. For all $a,b,c>1$, 
there exists $C_0,h_0>0$ such that for any $f\in L^2(A)$
$$
\eee_{a\bar{h},b\tilde{h},c h}^A(f)\leq C_0\ \eee_{\bar{h},\tilde{h},h}^A(f)
$$
for all $h\in]0,h_0]$.
\end{lem}
\bp
This is similar to the proof of Lemma 2.1 in \cite{DiLeMi11}. We leave it to the reader.
\ep

As for the Sobolev norm, we introduce the vectorial parameter $\vh=(\bar{h},\tilde{h},h)$ and we denote by 
$\vh \cdot x:=(\bar h x_1,\tilde h x', hx'')$ the inhomogenous action of $\vh$ on $x\in\R^d$.
We will also denote $\vh^{-1}=(\bar{h}^{-1},\tilde{h}^{-1},h^{-1})$, $\vh^d=\bar{h}\tilde{h}^{d'}h^{d''}$ and 
$\eee_\vh^A(g)=\eee_{\bar{h},\tilde{h},h}^A(g)$. With these notations, one has
\begin{equation*}
\begin{split}
\eee_\vh^A(g)&=\frac{1}{2\Vr_d \vh^{d}}
\int_{A\times A}\indic_{ \lvert x_1-y_1\lvert <\bar{h}, \lvert x'-y'\lvert <\tilde{h}, \lvert x''-y''\lvert <h}\lvert g(x)-g(y)\lvert ^2 dxdy\\
&=\frac{1}{2\Vr_d \vh^{d}}
\int_{A\times A}\indic_{|\vh^{-1}\cdot (x-y)|_\infty<1}|g(x)-g(y)|^2dxdy
\end{split}
\end{equation*}
\begin{lem}\label{lem:dirich-cube}
There exists some constants $C>1$ and $h_0>0$ such that for any $f\in L^2(\Omega_0)$, one has
$$
\frac 1C \eee^{ \Omega_0}_{\vh}(f)\leq \eee^{B_0}_{\vh}(f\circ\theta^{-1})\leq C \eee^{ \Omega_0}_{\vh}(f)
$$
for all $\vh=(\bar{h},\tilde{h},h)$ such that $0<h,\tilde{h}<h_0$ and $0<\bar{h}\leq\tilde {h}$.
\end{lem}
\bp
Since the jacobian of $\theta^{-1}$ is bounded one has
\begin{equation*}
\begin{split}
\eee_{\vh}^{ \Omega_0}(f)&=\frac{1}{2\Vr_d \vh^d}
\int_{ \Omega_0\times\Omega_0}\indic_{ \lvert x_1-y_1\lvert <\bar{h}, \lvert x'-y'\lvert <\tilde{h}, \lvert x''-y''\lvert <h}\lvert f(x)-f(y)\lvert ^2 dxdy\\
&\leq \frac{1}{2\Vr_d  \vh^d}
\int_{B_0\times B_0}\indic_{ \lvert x_1-y_1\lvert <\bar{h}, \lvert x_1^\alpha x'-y_1^\alpha y'\lvert <\tilde{h}, \lvert x''-y''\lvert <h}\\
&\phantom{*********************}\lvert f\circ\theta^{-1}(x)-f\circ\theta^{-1}(y)\lvert ^2 dxdy.
\end{split}
\end{equation*}
On the other hand, since $\alpha\geq 1$ then for $x,y\in B_0$ such that $\lvert x_1-y_1\lvert <\bar{h}\leq\tilde h$, one has 
$$
\lvert x_1^\alpha x'-y_1^\alpha y'\lvert\geq \vert x_1\vert^\alpha\vert x'-y'\vert-\lvert x_1^\alpha-y_1^\alpha\lvert \vert y'\vert\geq 2^{-\alpha}\vert x'-y'\vert-C_\alpha \tilde h
$$
for some constant $C_\alpha>0$.
This implies that 
\begin{equation*}
\begin{split}
\eee_{\vh}^{ \Omega_0}(f)
\leq \frac{1}{2\Vr_d  \vh^d}
\int_{B_0\times B_0}&\indic_{ \lvert x_1-y_1\lvert <\bar{h}, \lvert  x'- y'\lvert <M_{\alpha}\tilde{h}, \lvert x''-y''\lvert <h}\\
&\phantom{*****}
\lvert f\circ\theta^{-1}(x)-f\circ\theta^{-1}(y)\lvert ^2 dxdy
\end{split}
\end{equation*}
with $M_{\alpha}=2^\alpha(1+C_\alpha)$. Since $\partial B_0$ is Lipschitz, it follows from  Lemma \ref{lem:mult-param-dirich} that
$
\eee_{\vh}^{ \Omega_0}(f)\leq \tilde C_\alpha\eee_{\vh}^{B_0}(f\circ\theta^{-1}),
$
which proves the left inequality. The right inequality is proved similarly.
\ep

Since for any $k\geq 0$, one has $\bar h_k=2^k h\leq 2^{k\alpha}h=\tilde h_k$, it follows from  Lemma \ref{lem:dirich-cube} and \eqref{eq:decomp-dirich-scales2},  that 
\be\label{eq:decomp-dirich-scales2bis}
\ooo(h^2)=\eee_h(f)\geq \frac 1 C\sum_{k\geq0}2^{-k(1+\alpha d')}\eee^{B_0}_{\bar h_k,\tilde h_k,h}(f\circ \sigma_k^{-1}).
\ee
Let $Q:=]0,1[^d$ and define the change of variable 
\be\label{eq:defbeta}
\beta: B_0\rightarrow Q
\ee
 given by $\beta(x_1,x',x'')=(2x_1-1, x',x'')$.
Working as in Lemma \ref{lem:dirich-cube}, we show that there exists a constant $C>1$ such that 
\be\label{eq:chvar-dir}
\frac 1 C\eee^Q_{\vh}(g\circ\beta^{-1})\leq \eee^{B_0}_{\vh}(g)\leq C\eee^Q_{\vh}(g\circ\beta^{-1}).
\ee
Combined with  \eqref{eq:decomp-dirich-scales2bis} this implies  that there exists $C_0,h_0>0$ such that  for $0<h<h_0$
\be\label{eq:decomp-dirich-scales3}
\ooo(h^2)=\eee_{ h}(f)\geq \frac 1{C_0} \sum_{k\geq 0}2^{-k(1+\alpha d')}\eee^{Q}_{\bar{h}_k,\tilde{h}_k,h}(f\circ \sigma_k^{-1}\circ\beta^{-1}).
\ee

Let us now study the Dirichlet form on the cube $Q$. 
For any $i=1,\ldots,d$ let $s_i$ denote the symmetry  with respect to the hyperplane $\{x_i=1\}$ and let $G$ be the abelian group  generated by the $s_i$. The group $G$ acts on $]0,2[^d$ and for every function $f\in L^2(Q)$, one can then define $g\in L^2(]0,2]^d)$ by $g_{\lvert ]0,1[^d}=f$ and for all 
$s \in G$, $g\circ s = g$ (we do not specify the value of $g$ on the hyperplanes $\{x_i=1\}$ since they are negligible sets).
Eventually, this permits to extend the function $g$  (by means of translations) to a $(2\mathbb{Z})^d$-periodic function on $\mathbb{R}^d$.  We then denote
$$
\begin{array}{c}
E:L^2(Q)\rightarrow L^2(\mathbb{T}^d)\\
f\mapsto g
\end{array}
$$
where $\mathbb{T}^d=(\mathbb{R}/2\mathbb{Z})^d$. From the preceding discussion,  $E$ is continuous from $L^2(]0,1[^d)$ into $L^2(\mathbb{T}^d)$ and from $H^1(]0,1[^d)$ into $H^1(\mathbb{T}^d)$.
Given $0\leq a<b\leq 1$, we denote 
$$\Pi^{d-1}_{]a,b[}=]a,b[\times \Pi^{d-1}.$$
We can  perform a partial periodization by using only symmetries with respect to hyperplanes $\{x_i=1\}$ with $i\geq 2$.
We obtain an extension map 
\begin{equation*}
\begin{array}{c}
\tilde E_{]a,b[}:L^2(]a,b[\times ]0,1[^{d-1})\rightarrow L^2(\Pi^{d-1}_{]a,b[}).\\
\end{array}
\end{equation*}
We also  introduce the following restriction operators 
\begin{equation}\label{eq:defrestrict}
\begin{split}
&R:L^2(\Pi^d)\rightarrow L^2( ]0,1[^{d})\\
&R_{]a,b[}:L^2(\Pi^d)\rightarrow L^2(]a,b[\times ]0,1[^{d-1}),\\
&R^1_{]a,b[}:L^2(\Pi^d)\rightarrow L^2(\Pi^{d-1}_{ ]a,b[}),\\
&\tilde R_{ ]a,b[}:L^2(\Pi^{d-1}_{]a,b[})\rightarrow L^2(]a,b[\times ]0,1[^{d-1}),
\end{split}
\end{equation}
which satisfy the following relations:
\be\label{eq:restrictext}
RE=\Id,\;\tilde R_{ ]a,b[}\tilde E_{]a,b[}=\Id,\;R_{ ]a,b[}=\tilde R_{]a,b[}R^1_{]a,b[}, R^1_{]a,b[}E=\tilde E_{]a,b[}.
\ee
Eventually, we observe that all these operators are continuous on  $H^1$ and $L^2$ spaces.
In order to get rid of boundary problems, the general idea  is now to compare the Dirichlet form $\eee^Q_{\bar{h},\tilde{h},h}$ with a suitable Dirichlet form on the torus. We first introduce the Metropolis operator on $\Pi^d$, defined by 
\begin{equation}\label{eq:metrotorus}
 \bar{T}_{\vh}(g)(x)= \frac{1}{\Vr_{\infty,d}  \vh^d} \int_{\mathbb{T}^d} \indic_{ \lvert x_1-y_1\lvert <\bar{h}, \lvert x'-y'\lvert <\tilde{h}, \lvert x''-y''\lvert <h}g(y)dy
\end{equation}
for any $ g\in L^2(\mathbb{T}^d)$, 
where $\Vr_{\infty,d}=\int_{\mathbb{T}^d} \indic_{ \lvert y_1\lvert <1, \lvert y'\lvert <1, \lvert y''\lvert <1}dy$.
The associated Dirichlet form is
\begin{equation*}
\begin{split}
\bar\eee_{\vh}(g)&:= \<(1-\bar{T}_{\vh})(g),g\>_{L^2( \mathbb{T}^d)}\\
&=\frac{1}{2\Vr_{\infty,d}  \vh^d} 
\int_{\mathbb{T}^d\times \mathbb{T}^d}
\indic_{|\vh^{-1}\cdot(x-y)|_\infty <1}\lvert g(x)-g(y)\lvert ^2 dxdy.
\end{split}
\end{equation*}

\begin{lem}\label{lem:extension-dirich}
There exists $C,h_0>0$ such that for all $0<h,\bar h ,\tilde h<h_0$ and all $f\in L^2(Q)$
$$\eee^Q_{\vh}(f)\leq 
\bar\eee_{\vh}(E(f))
\leq C\eee^Q_{\vh}(f).$$
\end{lem}

\bp For any $f\in L^2(Q)$, one has
\begin{equation*}
\begin{split}
\bar\eee_{\vh}(E(f))&=\frac{1}{2\Vr_{\infty,d} \vh^d} 
\int_{\mathbb{T}^d\times \mathbb{T}^d}
\indic_{ \lvert \vh^{-1}\cdot(x-y)\vert_\infty<1}
\lvert g(x)-g(y)\lvert ^2 dxdy\\
& =
\frac{1}{2\Vr_{\infty,d}\vh^d}  \sum_{s,\tilde s\in G} \int_{s(Q)\times \tilde s(Q)}
\indic_{ \lvert \vh^{-1}\cdot(x-y)\vert_\infty<1}
\lvert g(x)-g(y)\lvert ^2 dxdy\\
& = 
\frac{1}{2\Vr_{\infty,d} \vh^d}  \sum_{s,\tilde s\in G} \int_{Q\times Q}
\indic_{  \lvert \vh^{-1}\cdot(s(x)-\tilde s(y))\vert_\infty<1}
\lvert g\circ s(x)-g\circ \tilde s(y)\lvert ^2
dxdy
\end{split}
\end{equation*}
and by definition of $g$ it follows 
\begin{equation}\label{eq:jesaispas}
\begin{split}
\bar\eee_{\vh}(E(f))
&=\frac{1}{2\Vr_{\infty,d} \vh^d}  \sum_{s,\tilde s\in G} \int_{Q\times Q}
\indic_{  \lvert \vh^{-1}\cdot(s(x)-\tilde s(y))\vert_\infty<1}
\\&\phantom{*****************}
\lvert f(x)-f(y)\lvert ^2
dxdy.
\end{split}
\end{equation}
Moreover,  for all $s,\tilde s\in G$  and any  $x,y\in Q$, one has
\be\label{eq:2star}
|\vh^{-1}\cdot(s(x)-\tilde s(y))\vert_\infty<1\Longrightarrow \lvert \vh^{-1}\cdot(x-y)\vert_\infty<2.
\ee
Indeed, since the elements of $G$ are isometries of $\R^d$ for the norm $|\cdot|_\infty$, it suffices to prove \eqref{eq:2star} with $s=\Id$.
If $\tilde s=\Id$, there is nothing to prove. Let us assume that $\tilde s\neq \Id$.
Then, there exists   $ a\in \{0,1\}^d$ such that $\tilde s=\prod_{i=1}^d s_i^{a_i}$.
Let us denote $I=\{i,\; a_i=1\}$ and let $D=\cap_{i\in I}\{x_i=1\}$. Since $x\in ]0,1[^d$, $\tilde s(y)\notin ]0,1[^d$ and $|\vh^{-1}\cdot(x-\tilde s(y))\vert_\infty<1$ then there exists $z\in D$ such that 
$|\vh^{-1}\cdot(x-z)\vert_\infty<1$ and $|\vh^{-1}\cdot(z-\tilde s(y))\vert_\infty<1$. Since $\tilde s(z)=z$,  this last inequality implies that 
$|\vh^{-1}\cdot(z-y)\vert_\infty<1$ and hence
$$
 \lvert \vh^{-1}\cdot(x-y)\vert_\infty\leq  \lvert \vh^{-1}\cdot(x-z)\vert_\infty+ \lvert \vh^{-1}\cdot(z-y)\vert_\infty<2
$$
which proves \eqref{eq:2star}.
Now, using  \eqref{eq:jesaispas} and \eqref{eq:2star}, we obtain
\begin{equation*}
\begin{split}
\bar\eee_{\vh}(E(f))
&\leq\frac{1}{2\Vr_{\infty,d} \vh^d}  \sum_{s,\tilde s\in G} \int_{Q\times Q}
\indic_{ \lvert \vh^{-1}\cdot(x-y)\vert_\infty<2}
\lvert f(x)-f(y)\lvert ^2dx dy
\end{split}
\end{equation*}
and thanks to Lemma \ref{lem:mult-param-dirich}, there exists $C,h_0>0$ such that for all $0<h<h_0$, one has 
$$
\bar\eee_{\vh}(E(f))\leq C \eee^Q_{\vh}(f).
$$
This proves the right inequality. The left one is immediate.
\ep

\begin{rem}\label{rem:extension-dirich}
The above proof can be easily adapted to show that given $0\leq a<b\leq 1$, there exists $C>0$ 
such that for all $0<h,\bar h ,\tilde h<h_0$ and all $f\in L^2(Q)$
$$
\bar\eee^{]a,b[}_{\vh}(\tilde E_{]a,b[}(f))
\leq C\eee^Q_{\vh}(f)$$
where 
\begin{equation*}
\begin{split}
\bar\eee^{]a,b[}_{\vh}(g)=\frac{1}{2\Vr_{\infty,d} \vh^d} 
\int_{\Pi_{]a,b[}^{d-1}\times\Pi_{]a,b[}^{d-1} }
&\indic_{ \lvert x_1-y_1\lvert <\bar{h}, \lvert x'-y'\lvert <\tilde{h}, \lvert x''-y''\lvert <h}\\
&\phantom{***********}
\lvert g(x)-g(y)\lvert ^2 dxdy.
\end{split}
\end{equation*}
\end{rem}
\subsubsection{Fourier analysis of the Metropolis operator on the torus}
The following lemma gives an expression of the operator $\bar{T}_{\vh}$ as a Fourier multiplier.
\begin{lem}\label{lem:Fourier-Th}
For $0<\bar h,\tilde h, h<1$, one has 
$$\bar{T}_{\vh}=\Gamma_1(\bar{h}^2\partial_1^2)\Gamma_{d'}(\tilde {h}^2\Delta_{x'})\Gamma_{d''}(h^2\Delta_{x''})$$
with 
$$\Gamma_n(\vert \xi\vert^2):=G_n(\xi):=\frac 1{\Vr_n}\int_{\R^n}\indic_{\vert z\vert<1}e^{i\pi z\cdot\xi}dz$$
\end{lem}
\bp 
First, observe that since $\Vr_{\infty,d}=\Vr_1\Vr_{d'}\Vr_{d''}$, one has 
$
\bar{T}_{\bar{h},\tilde{h}, h}=M_{1,\bar h}M_{d',\tilde h}M_{d'',h}
$
where for any $\hbar>0$ we set
$$
M_{n,\hbar}g(x)= \frac{1}{\Vr_{n} \hbar^n} \int_{\mathbb{T}^n} \indic_{ \lvert x-y\lvert <\hbar}g(y)dy.
$$
On the other hand, if one denotes $e_k:=\frac1{2^{n/2}}e^{i\pi z\cdot k}$ for all $k\in \Z^n$, then 
$(e_k)$ is an orthonormal basis of $L^2((\R/2\Z)^n)$. 
Moreover, for $0<\hbar<1$, the map $y\mapsto x+\hbar y$ is a change of variable on $\Pi^n$ and we get  
$$
M_{n,\hbar}(e_k)=e_k(x) \frac{1}{\Vr_{n} \hbar^n} \int_{\mathbb{T}^n} \indic_{ \lvert x-y\lvert <\hbar}e^{i\pi  k\cdot (y-x)}dy= G_n(\hbar k) e_k.
$$
Since the function $G_n$ is radial, this proves the announced result.
\ep

From the discussion below (1.6) in \cite{LeMi10} one knows that $G_n$ is a  smooth functions on $\R^n$ such that $\vert G_n\vert\leq 1$, $\vert G_n(\xi)\vert=1$ iff $\xi=0$ and 
\be\label{eq:DLGn}
G_n(\xi)=1-\frac 1{2(n+2)}\vert \xi\vert^2+\ooo(\vert \xi\vert^4).
\ee
With the notation \eqref{eq:notnormdyad}, we have the following
\begin{lem}\label{lem:regularite-tore} There exists $C, h_0>0$ such that for all  $0<h,\bar h,\tilde h<h_0$ and all
$g\in L^2(\Pi^d)$ such that $\Vert g\Vert_{L^2}\leq 1$ and $\bar \eee_{\vh}(g)\leq h^2$, there exists a decomposition 
$g=g_L+g_H$ such that
$$\Vert g_L\Vert_{L^2(\Pi^d)}^2+h^{-2}N_{\vh}(g_L,\Pi^d)  \leq C
$$
and
$$\lvert \lvert g_H\rvert \rvert_{L^2(\Pi^d)}^2\leq Ch^2.$$
\end{lem}

\begin{rem} Let $\lambda_1=\bar h/h$, $\lambda_2=\tilde h/h$.
One has
$$
\Vert g_L\Vert_{L^2(\Pi^d)}^2+h^{-2}N_{\vh}(g_L,\Pi^d)=\Vert g_L\Vert^2_{H^1_{\lambda_1,\lambda_2,1}(\Pi^d)}
$$
where the semiclassical Sobolev spaces $H^1_{\lambda_1,\lambda_2,1}$ are defined in appendix.
\end{rem}
\bp Denote $\alpha_n=\frac 1{2(n+2)}$ and let $1/\Upsilon_1=\frac 14\min(\alpha_1,\alpha_{d'},\alpha_{d''})>0$. 
Let $g\in L^2(\Pi^d)$  be such that  $\bar\eee_{\bar{h},\tilde{h},h}(g)\leq h^2$.
From  Lemma \ref{lem:Fourier-Th}, one knows that for any  $0<h,\bar h,\tilde h<1$, one has 
\be\label{eq:formDirG}
\begin{split}
h^2&\geq\<(1-\bar{T}_{\vh})g,g\>_{L^2(\Pi^d)}\\
&\geq\<(1-G_1(\bar{h}\partial_1)G_{d'}(\tilde{h}\nabla')G_{d''}(h\nabla''))g,g\>_{L^2(\Pi^d)}
\end{split}
\ee
On the other hand, it follows from \eqref{eq:DLGn} that  for all $n\in\N$, there exists $\delta_n>0$ such that for all $\vert\xi\vert<\delta_n$, one has
$$
0<G_n(\xi)\leq 1-\frac {\alpha_n}2\vert\xi\vert^2.
$$
Hence, for all $\xi=(\xi_1,\xi',\xi'')\in\R^d$ such that $\vert\xi\vert<\delta:=\min(\delta_1,\delta_{d'},\delta_{d''})$, one has
\begin{equation*}
\begin{split}
1-G_1(\xi_1)G_{d'}(\xi')G_{d''}(\xi'')&\geq \frac {\alpha_1}2\vert\xi_1\vert^2+\frac {\alpha_{d'}}2\vert\xi'\vert^2+\frac {\alpha_{d''}}2\vert\xi''\vert^2
+\ooo(\vert\xi\vert^4)\\
&\geq \frac 2{\Upsilon_1}\vert \xi\vert^2
+\ooo(\vert\xi\vert^4).
\end{split}
\end{equation*}
Decreasing $\delta$ as much as necessary, we obtain 
\be\label{minorG_norigin}
1-G_1(\xi_1)G_{d'}(\xi')G_{d''}(\xi'')\geq \frac 1{\Upsilon_1}\vert \xi\vert^2
\ee
for any $\vert\xi\vert<\delta$. On the other hand, since $G_n$  is bounded by $1$ and goes to zero at infinity and $1-G_n$ vanishes only at the origin, there exists $\Upsilon_2>0$ such that for all 
$\vert\xi\vert\geq \delta$, 
\be\label{minorG_ninfini}
1-G_1(\xi_1)G_{d'}(\xi')G_{d''}(\xi'')\geq \frac 1{\Upsilon_2}.
\ee
Let us decompose $g$ in the Fourier basis $(e_k)$,
$
g=\sum_{k\in\Z}\hat g(k) e_k
$
and let 
\begin{equation*}
g_L:=\sum_{\vert (\bar{h}k_1,\tilde{h}k',hk'')\vert<\delta}\hat g(k)e_k,\;\;g_H:=1-g_L=\sum_{\vert (\bar{h}k_1,\tilde{h}k',hk'')\vert\geq\delta}\hat g(k)e_k.
\end{equation*}
From \eqref{eq:formDirG}, \eqref{minorG_norigin} and \eqref{minorG_ninfini}, one deduces
\begin{equation*}
\begin{split}
h^2&\geq\sum_{k\in\Z^d}(1-G_1(\bar{h}k_1)G_{d'}(\tilde{h}k')G_{d''}(hk''))\vert\hat g(k)\vert^2\\
&=\sum_{\vert (\bar{h}k_1,\tilde{h}k',hk'')\vert<\delta}(1-G_1(\bar{h}k_1)G_{d'}(\tilde{h}k')G_{d''}(hk''))\vert\hat g(k)\vert^2\\
&\phantom{*********}+\sum_{\vert (\bar{h}k_1,\tilde{h}k',hk'')\vert\geq\delta}(1-G_1(\bar{h}k_1)G_{d'}(\tilde{h}k')G_{d''}(hk''))\vert\hat g(k)\vert^2\\
&\geq \frac 1{\Upsilon_1}\sum_{\vert (\bar{h}k_1,\tilde{h}k',hk'')\vert<\delta}(\vert \bar{h}k_1\vert^2+\vert \tilde{h}k'\vert^2+\vert hk''\vert^2)\vert\hat g(k)\vert^2\\
&\phantom{******************}+\frac 1{\Upsilon_2}\sum_{\vert (\bar{h}k_1,\tilde{h}k',hk'')\vert\geq\delta}\vert\hat g(k)\vert^2.
\end{split}
\end{equation*}
From standard Fourier analysis, we deduce
\begin{equation*}
h^2\geq \frac 1{\Upsilon_1}\Big(\Vert \bar{h}\partial_1g_L\Vert^2+\Vert\tilde{h}\nabla'g_L\Vert^2+\Vert h\nabla g_L\Vert^2\Big)
+\frac 1{\Upsilon_2}\Vert g_H\Vert^2
\end{equation*}
Taking $C=\max(\Upsilon_1,\Upsilon_2)$ we get the announced result.
\ep

\subsubsection{Decomposition Lemma in the cusp.}
The main result of this section is the following
\begin{lem}\label{lem:decompHFcusp}
Assume that $\Omega$ has the particular form \eqref{eq:defOmegamodel}, then the conclusion of Lemma \ref{lem:decompHF} holds true.
\end{lem}

\bp Throughout, $C$ denotes a positive constant independent of $f$ and $h$ that may change from line to line and $h>0$ is supposed sufficiently small in order that the conclusions of the preceding lemmas hold true.
Let $f\in L^2(\Omega)$ be such that
$\eee_h(f)\leq h^2\text{ and }
\lvert \lvert f\rvert \rvert _{L^2}\leq 1.$
It follows from \eqref{eq:decomp-dirich-scales3} that 
\begin{equation*}
\sum_{k\geq 0}2^{-k(1+\alpha d')}\eee^{Q}_{\bar{h}_k,\tilde{h}_k,h}(f\circ \sigma_k^{-1}\circ\beta^{-1})=\ooo(h^2)
\end{equation*}
where $\beta$ is defined by \eqref{eq:defbeta}.
Denoting 
\be\label{eq:defgbark}
\bar g_k=E(f\circ \sigma_k^{-1}\circ \beta^{-1}),
\ee
it follows from Lemma \ref{lem:extension-dirich} that
\be\label{eq:sum-dirich-part0}
\bar\eee_{\bar{h}_k,\tilde{h}_k,h}(\bar g_k)\leq C \eee^{Q}_{\bar{h}_k,\tilde{h}_k,h}(f\circ \sigma_k^{-1}\circ\beta^{-1})
\ee
and hence
\be\label{eq:sum-dirich-part}
\sum_{k\geq 0}2^{-k(1+\alpha d')}\bar\eee_{\bar{h}_k,\tilde{h}_k,h}(\bar g_k)=\ooo(h^2).
\ee
From now, given $D\subset\R^p\times \Pi^q$, $p,q\in\N^*$, $g\in L^2(D)$ and $\vh=(\bar h,\tilde h,h)$ we denote
\be\label{eq:defBrond}
\vvv^D_{\vh}(g)=\vvv^D_{\bar h,\tilde h,h}(g):=\Vert g\Vert_{L^2(D)}^2+h^{-2}\eee^D_{\vh}(g)
\ee
and for shortness we denote $\vvv_k(f,h)=\vvv^{B_0}_{\bar{h}_k, \tilde{h}_k,h}(f\circ\sigma_k^{-1})$.
We also denote $\vh_k=(\bar{h}_k, \tilde{h}_k,h)$.
Thanks to \eqref{eq:decomp-dyad-Sob1} and  \eqref{eq:sum-dirich-part}, one has 
\be\label{eq:estimBkf}
\sum_{k=0}^\infty 2^{-k(1+\alpha d')}\vvv_k(f,h)=\ooo(1),
\ee
and 
\eqref{eq:chvar-dir} and \eqref{eq:sum-dirich-part0} implies 
$$
\Vert \bar g_k\Vert^2_{L^2(\Pi^d)}+ h^{-2}\bar\eee_{\vh_k}(\bar g_k)\leq C\vvv_k(f,h).
$$
This estimate combined with  Lemma \ref{lem:regularite-tore} shows that there exists ${ h_0}>0$ such that for any $k\in\N$
 and $h>0$ such that $\tilde h_k <{ h_0}$, there exists $\bar g_{k,L}\in H^1(\Pi^d)$ and $\bar g_{k,H}\in L^2(\Pi^d)$ such that 
$\bar g_k=\bar g_{k,L}+\bar g_{k,H}$ with
$$
 \Vert \bar g_{k,L}\Vert_{L^2(\Pi^d)}^2+ h^{-2}N_{\vh_k}(\bar g_{k,L},\Pi^d)^2  \leq C  \vvv_k(f,h)
$$
and
$$
\Vert \bar g_{k,H}\Vert ^2_{L^2}\leq C h^2 \vvv_k(f,h)
$$
for some new constant $C$.
Since the restriction operator $R^1_{]0,1[}$ (defined in \eqref{eq:defrestrict}) is continuous, it follows from the above estimates that 
 \be\label{eq:estim-gkL}
 \Vert R^1_{]0,1[}(\bar g_{k,L})\Vert_{L^2( \Pi^{d-1}_{]0,1[} )}^2+ h^{-2}N_{\vh_k}(R^1_{]0,1[}(\bar g_{k,L}),\Pi^{d-1}_{]0,1[})^2  \leq C  \vvv_k(f,h)
 \ee
and
 \be\label{eq:estim-gkH}
 \lvert \lvert R^1_{]0,1[}(\bar g_{k,H})\rvert \rvert ^2_{L^2(\Pi^{d-1}_{]0,1[})}\leq C
 h^2 \vvv_k(f,h)
 \ee
 Combined with \eqref{eq:estimBkf}, this implies
\be\label{estimglobH1}
\sum_{k\geq 0}2^{-k(1+\alpha d')}\Vert R^1_{]0,1[}(\bar g_{k,H})\Vert ^2_{L^2(\Pi^{d-1}_{]0,1[})}\leq Ch^2
\ee
and
$$
\sum_{k\geq 0}2^{-k(1+\alpha d')}\Big(\Vert R^1_{]0,1[}(\bar  g_{k,L})\Vert_{L^2(\Pi^{d-1}_{]0,1[})}^2+h^{-2}N_{\vh_k} (R^1_{]0,1[}(\bar g_{k,L}),\Pi^{d-1}_{]0,1[})^2\Big)\leq C.
$$
Since $N_{\vh_k} (.,.)=hN_{2^k,2^{k\alpha},1}(.,.)$, this later equation implies
\be\label{estimglobL1}
\begin{split}
\sum_{k\geq 0}2^{-k(1+\alpha d')}\Big(\Vert R^1_{]0,1[}&(\bar  g_{k,L})\Vert_{L^2(\Pi^{d-1}_{]0,1[})}^2\\
&+N_{2^k,2^{k\alpha},1} (R^1_{]0,1[}( \bar g_{k,L}),\Pi^{d-1}_{]0,1[})^2\Big)=\ooo(1).
\end{split}
\ee
In view of Lemma \ref{lem:decomp-dyad-Sob2}, estimates \eqref{estimglobH1} and  \eqref{estimglobL1} almost imply the conclusion by considering the restriction of $R^1(\bar  g_k)$ to $B_0$. The main issue to get the conclusion is that nothing insures that the no-jump condition of the lemma
$R^1_{]0,1[}(\bar g_{k,L})_{\vert x_1=\frac 12}=R^1_{]0,1[}(\bar g_{k+1,L})_{\vert x_1=1}$ holds true (observe here that 
the interface $x_1=2^{-k-1}$ in the orginial variable corresponds to $x_1=\frac 12$ for $\bar g_{k,L}$ and to $x_1=1$ for $\bar g_{k+1,L}$).
The end of the proof consists to modify slightly the above decomposition in order to satisfy the assumptions of Lemma \ref{lem:decomp-dyad-Sob2}.

Let us explain briefly the idea of this modification before entering into technical details. As already said, 
it follows from the above estimates that we can decompose 
the functions $f_{|\Omega_{2k}}$ and $f_{|\Omega_{2k+1}}$  as the sum of a $H^1$ function and  a small function in $L^2$.
The idea is that we can do an analogous decomposition with quadriadic decomposition so that the function 
$f$ restricted to $\Omega_{2k}\cup\Omega_{2k+1}$ admits also a decomposition. Then we can apply Lemma \ref{lem:recol} of the appendix to 
glue smoothly $f_{|\Omega_{2k}}^L$ and $f_{|\Omega_{2k+1}}^L$ up to a small error in $L^2$. In order to get a global estimate, 
we need to prove estimates uniform with respect to the dyadic parameter $k$ which makes the computation a bit more heavy.

Let us now enter into the details. We first observe that thanks to \eqref{eq:restrictext}, one has
 \be\label{eq:splitf1}
 f
 =\sum_{k\geq 0}\indic_{\Omega_k} f_k\circ\sigma_k=\sum_{k\geq 0}\indic_{ \Omega_k} R(\bar g_k)\circ\beta\circ\sigma_k
\ee
 with $f_k=f\circ\sigma_k^{-1}=R(\bar g_k)\circ\beta$ ($\bar g_k$ given by \eqref{eq:defgbark}, $\beta$ given by \eqref{eq:defbeta} and $R$ given by \eqref{eq:defrestrict}).
 We introduce the following functions defined on   $\Pi^{d-1}_{]\frac 12,1[}$: 
\be\label{eq:def-gcheck}
\check g_k=R^1_{]0,1[}(\bar g_k)\circ\beta,\;\check g_{k,L}=R^1_{]0,1[}(\bar g_{k,L})\circ\beta,\;
\check g_{k,H}=R^1_{]0,1[}(\bar g_{k,H})\circ\beta
\ee
which of course verify
\be\label{eq:decompgcheck}
\check g_K=\check g_{k,L}+\check g_{k,H}
\ee
thanks to the above construction.
First observe that thanks to \eqref{eq:restrictext}, one has for any $k\in\N$, 
$f_k=R(\bar g_k)\circ\beta=\tilde R_{]0,1[}(R^1_{]0,1[}(\bar g_k))\circ\beta=\tilde R_{]\frac 12,1[}(R^1_{]0,1[}(\bar g_k)\circ\beta)$ and hence
\be\label{eq:pulfk1}
f_k=\tilde R_{]\frac 12,1[}(\check g_k)
\ee
where we recall that $\tilde R_{]a,b[}:L^2(\Pi_{]a,b[}^{d-1})\rightarrow L^2(]a,b[\times]0,1[^{d-1})$ denotes the restriction operator in the $(x',x'')$ variable.
Splitting  \eqref{eq:splitf1} into even and odd terms and using \eqref{eq:pulfk1}, we get
\be\label{eq:splitf2}
\begin{split}
f&=\sum_{k\geq 0}\indic_{ \Omega_{2k}}\tilde R_{]\frac 12 ,1 [}(\check g_{2k})\circ\sigma_{2k}\\
&\phantom{***********}
+\sum_{k\geq 0}\indic_{ \Omega_{2k+1}}\tilde R_{]\frac 12 ,1[}(\check g_{2k+1})\circ \check\sigma_1\circ\check\sigma_1^{-1}\circ\sigma_{2k+1}
\end{split}
\ee
with $\check\sigma_1$ given by 
\eqref{eq:defchecksigma1}. Since 
this change of variable is simply given by
 $\check\sigma_1(x)=(2x_1,x',x'')$ we have for any $\psi$
$$\tilde R_{]\frac 12,1[}(\psi)\circ\check\sigma_1
=\tilde R_{]\frac 14 ,\frac 12[} (\psi\circ\check\sigma_1)$$
where with a slight abuse of notation we use the symbol $\check\sigma_1$ to denote the above dilation defined from 
$\Pi^{d-1}_{]\frac 14,\frac 12[}$ into $\Pi^{d-1}_{]\frac 12,1[}$.
Combined with \eqref{eq:splitf2} and the identity 
$\check\sigma_1^{-1}\circ\sigma_{2k+1}=\hat\sigma_{2k+1}$ (see \eqref{eq:def-hattauk}, \eqref{eq:def-hatsigmak} for the definition of $\hat\sigma_k$), this implies
 \be\label{eq:splitf3}
\begin{split}
f=\sum_{k\geq 0}\indic_{\Omega_{2k}}&\tilde R_{]\frac 12 ,1[}(\check g_{2k})\circ\sigma_{2k}\\
&+\sum_{k\geq 0}\indic_{ \Omega_{2k+1}}\tilde R_{]\frac 14 ,\frac 12[}(\check g_{2k+1}\circ \check\sigma_1)\circ\hat \sigma_{2k+1}.
\end{split}
\ee
Denote $D_k= \Omega_{2k}\cup \Omega_{2k+1}$ for any $k\in \N$ and let 
$ \nu_k:\;\, D_k\rightarrow B_0\cup B_1$ be defined by 
 $\nu_k(x)=\theta(4^k x_1,4^{\alpha k}x',x'')$.
Since $(\nu_k)_{\vert\Omega_{2k}}=\sigma_{2k}$ and $(\nu_k)_{\vert\Omega_{2k+1}}=\hat\sigma_{2k+1}$, equation  \eqref{eq:splitf3} becomes
 \be\label{eq:splitf4}
\begin{split}
f=\sum_{k\geq 0}\indic_{ \Omega_{2k}}&\tilde R_{]\frac 12 ,1[}(\check g_{2k})\circ\nu_{k}\\
&+\sum_{k\geq 0}\indic_{\Omega_{2k+1}}\tilde R_{]\frac 14 ,\frac 12[}(\check g_{2k+1}\circ \check\sigma_1)\circ\nu_k.
\end{split}
\ee
We now relate the quadriadic decomposition of the function $f$ to the dyadic decomposition.
\begin{slem}\label{lem:quadriadyc}
Let $\check f_k:=\tilde E_{]\frac 14, 1[}(f\circ\nu_k^{-1})$ and denote 
\be\label{eq:defVrk}
\tilde \vvv_k(f,h)=\vvv_{\bar h_{2k}\tilde h_{2k},h}^{B_0\cup B_1}(f\circ\nu_k^{-1})
\ee
where the functional  $\vvv$ is defined by \eqref{eq:defBrond}.
One has 
\be\label{eq:retourvarinit4}
\begin{split}
\check f_k&=\indic_{\frac 14< x_1<\frac 12}\check g_{2k+1}\circ\check\sigma_1
+\indic_{\frac 12< x_1< 1}\check g_{2k}
\end{split}
 \ee
and 
\be\label{eq:estimBkfquadri}
\sum_{k=0}^\infty 4^{-k(1+\alpha d')}  \tilde\vvv_k(f,h)=O(1).
\ee
Moreover, one has the decompositions $\check g_{\bullet}=\check g_{\bullet,L}+\check g_{\bullet,H}$
and $\check f_k=\check f_{k,L}+\check f_{k,H}$ with 
\begin{equation}\label{eq:estimcheckg2k}
\begin{split}
&\Vert\check g_{2k+1,L}\circ\check\sigma_1\Vert^2_{L^2(\Pi_{]\frac 14,\frac 12[}^{d-1})}+ h^{-2}N_{\vh_{2k}} (\check g_{2k+1,L}\circ\check\sigma_1,\Pi_{]\frac 14,\frac 12[}^{d-1} )^2\\
&\phantom{*************************}\leq C 
\tilde \vvv_k(f,h)\\
&\Vert \check g_{2k,L}\Vert^2_{L^2(\Pi_{]\frac 12,1[}^{d-1})}+h^{-2}N_{\vh_{2k}} ( \check g_{2k,L},\Pi_{]\frac 12,1[}^{d-1} )^2\leq C
\tilde \vvv_k(f,h)
\end{split}
\end{equation}
and
\be\label{eq:estimNormD0H}
\Vert \check g_{2k+1,H}\circ\check\sigma_1\Vert ^2_{L^2(\Pi_{]\frac 14,\frac 12[}^{d-1})}
+\Vert \check g_{2k,H}\Vert ^2_{L^2(\Pi_{]\frac 12,1[}^{d-1})}
\leq C h^2 \tilde \vvv_k(f,h)
\ee
and
\be\label{eq:estimNormD2L}
\Vert \check f_{k,L}\Vert^2_{L^2(\Pi_{]\frac 14,1[}^{d-1})}
+ h^{-2}N_{\vh_{2k}} (\check f_{k,L},\Pi_{]\frac 14,1[}^{d-1} )^2
\leq  C
\tilde \vvv_k(f,h)
\ee
and
\be\label{eq:estimNormD2H}
\Vert \check f_{k,H}\Vert ^2_{L^2(\Pi_{]\frac 14,1[}^{d-1})}
\leq  Ch^2 \tilde \vvv_k(f,h).
\ee
where $C$ is a positive constant independent of $h$.
\end{slem}
\bp
By definition, one has 
$
f=\sum_{k\geq 0}\indic_{D_k} f\circ\nu_k^{-1}\circ\nu_k
$,
which combined to \eqref{eq:splitf4}  proves that 
\be\label{eq:splitf5}
f\circ\nu_k^{-1}=\indic_{B_0}\tilde R_{]\frac 12 ,1[}(\check g_{2k})+\indic_{B_1}\tilde R_{]\frac 14 ,\frac 12[}(\check g_{2k+1}\circ \check\sigma_1).
\ee
Applying $\tilde E_{[\frac 14, 1]}$ on both sides of this identity, we get \eqref{eq:retourvarinit4}.

We now observe that the analysis of Lemma \ref{lem:extension-dirich},\ref{lem:Fourier-Th}, \ref{lem:regularite-tore} can be performed with the "quadriadic" decomposition of the cusp induced by the change of variable $\nu_k$. This yields 
\be\label{eq:decomp-new-scales1}
h^2\geq\eee_h(f)\geq \frac 1C\sum_{k\geq0}4^{-k(1+\alpha d')}\eee^{B_0\cup B_1}_{4^k h,4^{k\alpha} h,h}(f\circ \nu_k^{-1}).
\ee
Dividing by $h^2$ and adding the $L^2$ norm, this implies
$$
\sum_{k=0}^\infty 4^{-k(1+\alpha d')}  \tilde\vvv_k(f,h)=O(1)
$$
which is exactly \eqref{eq:estimBkfquadri}.
Moreover,  it follows from Remark \ref{rem:extension-dirich}, \eqref{eq:pulfk1}, \eqref{eq:splitf5} and the inclusion 
$B_{0}\times B_{0}\cup B_{1}\times B_{1}\subset\check (B_0\cup B_1)^2$ that 
\be\label{eq:splitform1}
\begin{split}
\eee^{B_0\cup B_1}_{4^k h,4^{k\alpha} h,h}(f\circ \nu_k^{-1})&\geq 
\eee^{B_0}_{\bar h_{2k},\tilde h_{2k},h}(f\circ \sigma_{2k}^{-1})
+\eee^{B_1}_{\bar h_{2k},\tilde h_{2k},h}(f\circ \hat \sigma_{2k+1}^{-1})\\
&\geq 
\frac 1C\Big( \bar \eee^{]\frac 12,1[}_{\bar h_{2k},\tilde h_{2k},h}(\check g_{2k})
+ \bar \eee^{]\frac 14,\frac 12[}_{\bar h_{2k},\tilde h_{2k},h}(\check g_{2k+1}\circ\check\sigma_1)\Big)
\end{split}
\ee
for some constant $C>0$.
Observe that $\bar{h}_{2k+1}=2\bar{h}_{2k}, \tilde{h}_{2k+1}=2^\alpha\tilde{h}_{2k}$. 
Hence \eqref{eq:splitform1} proves  that there exists $C>0$ such that
\be\label{eq:estimVr}
\vvv_{2k}(f,h)+\vvv_{2k+1}(f,h)\leq C \tilde\vvv_k(f,h).
\ee 
Combined with   \eqref{eq:estim-gkL},  \eqref{eq:estim-gkH}, \eqref{eq:def-gcheck}, this proves 
\eqref{eq:estimcheckg2k} and 
\eqref{eq:estimNormD0H}.
On the other hand,  using Lemma \ref{lem:regularite-tore} and \eqref{eq:decomp-new-scales1} we have also a decomposition
$$E(f\circ\nu_k^{-1})=E(f\circ\nu_k^{-1})_L+E(f\circ\nu_k^{-1})_H$$
with suitable bounds on the right hand side. Restricting this decomposition to $\frac 14\leq x_1\leq 1$, we get 
$
\check f_k=\check f_{k,L}+\check f_{k,H}
$
with $\check f_{k,L}, \check f_{k,H}$ which satisfy 
\eqref{eq:estimNormD2L}
and
\eqref{eq:estimNormD2H}. This completes the proof of the sub-lemma.
\ep

Let us  now apply Lemma \ref{lem:recol} with $A_0=\Pi_{]\frac 14,\frac 12[}^{d-1}$, $A_1=\Pi_{]\frac 12,1[}^{d-1}$, $A_2=\Pi_{]\frac 14,1[}^{d-1}$, 
$\phi_0=\check g_{2k+1,L}\circ\check\sigma_1$ and $\phi_1=\check g_{2k,L}$.
Let $w_k:=\indic_{A_0}\phi_0+\indic_{A_1}\phi_1$ and denote  $r_0=\indic_{A_0} \check g_{2k+1,H}\circ\check\sigma_1$, 
$r_1=\indic_{A_1}\check g_{2k,H}$.
Thanks to \eqref{eq:decompgcheck} and  \eqref{eq:retourvarinit4}, $w_k$ satisfies
\be
\begin{split}
w_k&=\indic_{A_0} \check g_{2k+1}\circ\check\sigma_1+\indic_{A_1}\check g_{2k}
-r_0-r_1=\indic_{A_2}\check f_k-r_0-r_1.
\end{split}
\ee
Hence,
$w_k=\phi_2+r_2$ with 
\be 
\phi_2=\indic_{A_2}\check f_{k,L}\text{ and }r_2=\indic_{A_2}\check f_{k,H}-r_0-r_1.
\ee
Moreover, thanks to \eqref{eq:estimNormD0H}, \eqref{eq:estimNormD2L} and \eqref{eq:estimNormD2H}, one has
\be\label{eq:estim-wkH}
 \lvert \lvert  r_2\rvert \rvert ^2_{L^2(A_2)}\leq C h^2
\tilde \vvv_k
 \ee
and 
 $$
 N_{4^k,4^{k\alpha},1}(\phi_2,A_2)\leq C\tilde \vvv_k.
 $$
 where we write for shortness $\tilde \vvv_k=\tilde \vvv_k(f,h)$.
From Lemma \ref{lem:recol} with $\lambda_1=2^{2k}$, $\lambda'=2^{2k\alpha}$ and $\lambda''=1$, there exists $\Upsilon_1>0$ and 
${h_1}>0$ such that 
 for any $k$ such that $\bar h_{2k}\leq {h_1}$ (that is $2^{-2k}>h/{h_1}$), there exists  a function $\psi_{2k}$ supported in
 $\Pi_{]\frac 14,\frac 12[}^{d-1}\cap\{\frac 12\leq x_1\leq \frac 12+\bar h_{2k}\}$ such that $(\psi_{2k})_{\vert x_1=\frac 12}=(\phi_0)_{\vert x_1=\frac 12}
 -(\phi_1)_{\vert x_1=\frac 12}$ and 
 \begin{equation*}
 N_{4^k,4^{k\alpha},1} (\psi_{2k},A_1)^2 \leq \Upsilon_1\tilde\vvv_k
  \end{equation*}
and
 \begin{equation*}
 \Vert \psi_{2k}\Vert ^2_{L^2(A_1)}\leq \Upsilon_1 h^2 \tilde\vvv_k.
  \end{equation*}
  From now, we suppose that $\tilde h_{2k}<{ h_2}:=\min({h_0},{ h_1})$ with ${ h_0}$ given by Lemma \ref{lem:regularite-tore} and ${ h_1}$ by Lemma \ref{lem:recol}.
 We then rewrite $\check g_{2k}$ as $\check g_{2k}=\check g_{2k,L}^{mod}+\check g_{2k,H}^{mod}$ with 
 $\check g_{2k,L}^{mod}=\check g_{2k,L}+\psi_{2k}$ and $\check g_{2k,H}^{mod}=\check g_{2k,H}-\psi_{2k}$. By construction, we have
 \begin{align}\label{eq:contrace}
 \left\{\begin{array}{c}
 N_{4^k,4^{k\alpha},1} (\check g_{2k,L}^{mod},\Pi^{d-1}_{]\frac 12,1[})^2 \leq C \tilde\vvv_k\text{ and }\\
 \Vert \check g_{2k,H}^{mod} \Vert ^2_{L^2(\Pi^{d-1}_{]\frac 12,1[})}\leq Ch^2 \tilde\vvv_k\phantom{++++}\\
 (\check g_{2k,L}^{mod})_{x_1=\frac 12}=(\check g_{2k+1,L}\circ\check\sigma_1)_{x_1=\frac 12}.\phantom{+}
 \end{array}
 \right.
\end{align}

\noindent
Let $K(h)\in\N$ be the largest integer such that $4^{\alpha K(h)}\leq { h_2}/h$. Then for $k\leq K(h)$ the functions $\check g_{2k,L}^{mod}$ and $\check g_{2k,H}^{mod}$ are well-defined and we can introduce the decomposition $f:=f_L+f_H+f_\ccc$ with 
\begin{equation*}
\begin{split}
f_L= \sum_{k=0}^{K(h)}\Big(\indic_{ \Omega_{2k}}\tilde R_{]\frac 12,1[}(\check g_{2k,L}^{mod})\circ\sigma_{2k}+\indic_{ \Omega_{2k+1}}\tilde R_{]\frac 12,1[}(\check g_{2k+1,L})\circ\sigma_{2k+1}\Big),
\end{split}
\end{equation*}
\begin{equation*}
\begin{split}
f_H=\sum_{k=0}^{K(h)}\Big(\indic_{ \Omega_{2k}}\tilde R_{]\frac 12,1[}(\check g_{2k,H}^{mod})\circ\sigma_{2k}
+\indic_{\Omega_{2k+1}}\tilde R_{]\frac 12,1[}(\check g_{2k+1,H})\circ\sigma_{2k+1}\Big),
\end{split}
\end{equation*}
and 
\begin{equation*}\label{eq:defFC}
f_\ccc=f-f_L-f_H.
\end{equation*}
It follows from \eqref{eq:splitf4} and the definition of $K(h)$ that $f_\ccc$ is supported in $\{0<x_1< (h/ h_2)^{1/\alpha}\}$ which is the required property 
on $f_\ccc$.
On the other hand, we deduce from \eqref{eq:estimcheckg2k} and  \eqref{eq:contrace} that
\begin{equation*}
\begin{split}
\sum_{k=0}^{\infty}2^{-k(\alpha d'+1)}&N_{2^k,2^{k\alpha},1}(f_L\circ\sigma_k^{-1}\circ\theta^{-1},B_0)^2\\
&\leq C \sum_{k=0}^{K(h)}4^{-k(\alpha d'+1)}N_{4^k,4^{k\alpha},1} (\check g_{2k,L}^{mod},\Pi^{d-1}_{]\frac 12,1[})^2\\
&\phantom{****}+C \sum_{k=0}^{K(h)}4^{-k(\alpha d'+1)}N_{4^k,4^{k\alpha},1} (\check g_{2k+1,L}\circ\check\sigma_1,\Pi^{d-1}_{]\frac 14,\frac 12[})^2\\
&\leq C\sum_{k=0}^\infty 4^{-k(\alpha d'+1)}\tilde \vvv_k\leq C'
\end{split}
\end{equation*}
where the last inequality follows from \eqref{eq:estimBkfquadri} and $C'$ is a positive constant.
Hence, $f_L$ satisfies
\be\label{eq:normH1final}
\sum_{k=0}^{\infty}2^{-k(\alpha d'+1)}N_{2^k,2^{k\alpha},1}(f_L\circ\sigma_k^{-1},\Omega_0)^2=\ooo(1)
\ee
and thanks to \eqref{eq:contrace} the functions 
$\indic_{\Omega_{2k+1}} f_L$ and $\indic_{ \Omega_{2k}} f_L$ have the same trace on $x_1=2^{-2k-1}$. 
Working similarly near $x_1=2^{-2k}$, we can modify $(f_L)_{ \Omega_{2k}}$ in order that 
$\indic_{ \Omega_{2k}} f_L$ and $\indic_{\Omega_{2k-1}} f_L$  have the same trace on $x_1=2^{-2k}$.
Moreover, this new modification  is supported 
in $2^{-2k}[1-h_0,1]$. Hence, for $h_0>0$ small enough, it doesn't intersect the support of the modification $\psi_{2k}$ which is contained in $2^{-2k}[\frac 12,\frac 12+h_0]$.
Eventually, we modify also the function $\check g_{2K(h)+1}$ in order that $(f_L)_{|x_1=  4^{-K(h)-1}}=0$.
Consequently, the fonction $f_L$ that we obtain satisfies the assumptions of Lemma \ref{lem:decomp-dyad-Sob2} and it follows that 
$f_L\in H^1(\Omega)$ and $\Vert f_L\Vert_{H^1(\Omega)}=\ooo(1)$. 
The fact that $\Vert f_H\Vert_{L^2(\Omega)}=\ooo( h)$ follows immediately 
from \eqref{eq:estimBkfquadri}, \eqref{eq:estimNormD0H}, \eqref{eq:contrace} and Lemma \ref{lem:decomp-dyad-Sob2}.
\ep

\subsection{The general case}\label{sec:generalcase}
Suppose that $(f_h)_{h\in]0,1]}$ is a family of functions in $L^2(\Omega)$ such that $\Vert f_h\Vert_{L^2}=1$ and 
$\eee_h(f_h)=\ooo(h^2)$.
Let $J=I_c\cup I_r\cup\{0\}$ and for all $j\in J$, let $\ooo_j=\omega_j\cap\Omega$ where the 
$\omega_j,\;j\in I_c\cup I_r$ are defined in Assumption \ref{hyp1} and $\omega_0$ is a relatively compact open subset of $\Omega$ such that 
$\Omega\subset \cup_{j\in J}\ooo_j$.
Since $J$ is finite (independent of $h$), there exists $C>0$ such that  for any $f\in L^2$ one has 
\begin{equation}\label{eq:minorDirichcartes}
\begin{split}
\eee_h(f)\geq\frac 1C\sum_{i\in J}
 \eee_h^{\ooo_i}(f)
\end{split}
\end{equation}
with
$$
 \eee_h^{\ooo_i}(f):=\frac 1{2 h^{d}}\iint_{\ooo_i\times\ooo_i}\indic_{\vert x-y\vert <h}\vert f(x)-f(y)\vert^2d\mu_\rho^2(x,y).
$$
Let $(\chi_i)_{i\in J}$ be a family of non negative smooth functions such that $\supp(\chi_i)\subset\omega_i$ for all 
$i\in J$ and $\sum_{i\in J}\chi_i=1$ near $\Omega$.
For all $i\in J$, denote $f_{i,h}=(f_h)_{\vert\ooo_i}$ and observe that 
$$
f_h=\sum_{i\in J}\chi_if_{i,h}.
$$
Moreover, since $\eee_h(f_h)=\ooo(h^2)$ it follows from \eqref{eq:minorDirichcartes} that for all $i\in J$, $ \eee_h^{\ooo_i}(f_{i,h})=\ooo(h^2)$.
Suppose first that $i\in I_c$. In a suitable coordinate system, $\ooo_i$ has the form \eqref{eq:formcusp} and we can apply Lemma 
\ref{lem:decompHFcusp} to get the decomposition
\be\label{eq:decomposlocal}
f_{i,h}=\varphi_{i,h}+g_{i,h}+r_{i,h}
\ee
with 
$\Vert g_{i,h}\Vert_{L^2(\ooo_i)}=\ooo(h)$, $\supp(r_{i,h})\subset\Gamma_h$ (where $\Gamma_h$ is defined in Lemma \ref{lem:decompHF})   and $(\varphi_{i,h})_{h\in]0,1]}$ bounded in $H^1(\ooo_i)$.
On the other hand, it follows from Lemma 2.2 in \cite{DiLeMi11} that for any $i\in I_r\cup\{0\}$, \eqref{eq:decomposlocal} holds true with $r_{i,h}=0$.
As a consequence, we get a global decomposition  $f_h=\varphi_h+g_h+r_h$ with 
$$
\varphi_h=\sum_{i\in J}\chi_i\varphi_{i,h},\;\;g_h=\sum_{i\in J}\chi_ig_{i,h},\;\;r_h=\sum_{i\in I_c}\chi_ir_{i,h}.
$$
The functions $g_h$ and $r_h$ satisfy trivially the required properties. Since $\chi_i$ is supported in $\ooo_i$, one has the identity
$$
\nabla(\chi_i\varphi_{i,h})=\varphi_{i,h}\nabla\chi_i+\chi_i\nabla\varphi_{i,h}
$$
which permits easily to show that  $(\varphi_h)$ is bounded in $H^1$.
\section{Spectral analysis.}\label{sec4}

\subsection{Weak convergence of Dirichlet forms.}
We start this section with a lemma giving estimates of the Dirichlet form $\eee_h$ on $H^1$ fonctions.
Given a subset $U$ of $\Omega$, we use the notation 
$$
 \eee_h^{U}(u):=\frac 1{2 h^{d}}\iint_{U\times U}\indic_{\vert x-y\vert <h}|u(x)-u(y)|^2d\mu_\rho^2(x,y).
$$
\begin{lem}\label{lem:estimEhUH1} Suppose that the domain $\Omega$ satisfies Assumption \ref{hyp1}.
There exists $C>0$ and $h_0>0$ such that for any subset $U\subset\Omega$ and any $u\in H^1(\Omega)$, one has for all $h\in]0,h_0]$
$$
\eee_h^U(u)\leq Ch^2\Vert \nabla u\Vert^2_{L^2(U+B(0,Ch))}.
$$
\end{lem}
\bp
From Theorem 2, p 27 in \cite{MaPo97}, we know that  $\ccc^\infty(\Omega)\cap H^1(\Omega)$ is dense in $H^1(\Omega)$ for any open set $\Omega$.
Hence, we can assume that $u\in C^\infty(\Omega)$. 
Let $(\omega_i)_{i\in J}$, $J=I_c\cup I_r\cup\{0\}$ be a covering of $\Omega$ as in section \ref{sec:generalcase}.
For any $i\in J$, we denote $\omega_i^h=\omega_i+B(0,h)$. We have
\begin{equation}\label{eq:estimEhUH1-1}
\begin{split}
 \eee_h^{U}(u)&:=\sum_{j\in J}\frac 1{2 h^{d}}\iint_{U\cap\omega_j\times U}\indic_{\vert x-y\vert <h}|u(x)-u(y)|^2d\mu_\rho^2(x,y)\\
& =\sum_{j\in J}\frac 1{2 h^{d}}\iint_{U\cap\omega_j\times U\cap\omega_j^h}\indic_{\vert x-y\vert <h}|u(x)-u(y)|^2d\mu_\rho^2(x,y)\\
&\leq \sum_{j\in J}\eee_h^{U\cap\omega_j^h}(u)
 \end{split}
 \end{equation}
 For any $j\in J$, since $\rho$ is bounded, one has 
 $$
 \eee_h^{U\cap\omega_j^h}(u)\leq C h^{-d}\iint_{x,y\in U\cap \omega_j^{h}}\indic_{\vert x-y\vert<h}\vert u(x)-u(y)\vert^2dxdy.
 $$
 and using the change of variable $y=x+hz$ this implies
$$
 \eee_h^{U\cap\omega_j^h}(u)\leq C\int_{x\in U\cap \omega_j^{h}}\int_{z\in D_h}\vert u(x)-u(x+hz)\vert^2dxdz
$$
where $D_h=\{z\in B(0,1),\;x+hz\in U\cap\omega_j^h\}$.  Using local coordinates in $\omega_j^h$ one sees that there exists a piecewise smooth path  $ \gamma_{x,z,h}:[0,1]\rightarrow \Omega$ 
joining $x$ to $x+hz$ in $U\cap\omega_j^h$ such that in local coordinates $\gamma_{x,z,h}$ is the union of two straight lines
from $x=(x_1,x')$ to $(x_1+hz_1,x')$ and from $(x_1+hz_1,x')$ to $(x_1+hz_1,x'+hz')$.  In particular
there exist $C>0$ independent of $x,z,h$ such that  $\vert \dot \gamma_{x,z,h}(t)\vert\leq C h$ for all 
$t\in[0,1]$  and $d_x\gamma_{x,z,h}(t)=Id+O(h)$ uniformly with respect to $z$ and $t$. 
Hence, for any $t\in[0,1]$, $z\in B(0,1)$ and $h>0$ small enough, the map $\kappa_{t,z,h}:x\mapsto \gamma_{x,z,h}(t)$ is a change of variable from $\{x\in U\cap\omega_j^h,\; x+thz\in U\cap\omega_j^h\}$ onto a subset $V_j^h$ of $U\cap\omega_j^h+B(0,h)$.
By the fundamental theorem of analysis, it follows that
$$
 \eee_h^{U\cap\omega_j^h}(u)\leq C\int_{x\in U\cap \omega_j^{h}}\int_{z\in D_h}\vert \int_0^1 \dot\gamma_{x,z,h}(t)\cdot\nabla u(\gamma_{x,z,h}(t)dt\vert^2dxdz
$$
and thanks to the bound  $\vert \dot \gamma_{x,z,h}(t)\vert\leq C h$ we get
$$
 \eee_h^{U\cap\omega_j^h}(u)\leq C h^2\int_{x\in U\cap \omega_j^{h}}\int_{z\in D_h}\int_0^1\vert \nabla u(\gamma_{x,z,h}(t))\vert^2dtdzdx
$$
and using the change of variable $y=\kappa_{t,z,h}(x)$ it follows that
\begin{equation*}
\begin{split}
 \eee_h^{U\cap\omega_j^h}(u)&\leq C h^2\int_0^1\int_{z\in B(0,1)}\int_{y\in V_j^{h}}\vert \nabla u(y)\vert^2dydzdt\\
 &\leq C h^2\int_{y\in U\cap\omega_j^h+B(0,h)}\vert \nabla u(y)\vert^2dy.
\end{split}
\end{equation*}
Plugging this inequality in the last inequality of \eqref{eq:estimEhUH1-1} and since $J$ is finite, we get the result.
\ep
From now, given $r>0$, we denote 
\be\label{eq:defOmegar}
\Omega_r=\{x\in\Omega,\;d(x,\partial\Omega)\leq r\}.
\ee
Since we do not use the notation $\Omega_k$ related to the dyadic decomposition of the cusp, there is no ambiguity.
\begin{cor}\label{cor:locEh}
Suppose that $(u_h)$ is a family of functions which is bounded in $H^1(\Omega)$. Then $\eee_h(u_h)=O(h^2)$.
Moreover, if $u\in H^1(\Omega)$ is a fixed function independent of $h$, one has 
$$
\eee_h(u)=\eee_h^{\Omega\setminus\Omega_h}(u)+o(h^2).
$$
\end{cor}
\bp
The first estimate is a direct consequence of the preceding lemma with $U=\Omega$.
To get the second estimate observe that 
$$
\eee_h(u)=\eee_h^{\Omega\setminus\Omega_h}(u)+\eee_h^{\Omega_h}(u)+R_h
$$
with 
$$
R_h=\frac 1{ h^{d}}\iint_{\Omega_h\times \Omega\setminus\Omega_h}\indic_{\vert x-y\vert <h}|u(x)-u(y)|^2d\mu_\rho^2(x,y)
\leq 2\eee_h^{\Omega_{2h}}(u).
$$
From the preceding lemma  with $U=\Omega_{2h}$, it follows that
$$
h^{-2}\eee_h^{\Omega_{2h}}(u)\leq C\int \indic_{\Omega_{3h}}(x)|\nabla u(x)|^2dx
$$
which goes to $0$ as $h\rightarrow 0$ by the dominated convergence theorem (since $u$ doesn't depend on $h$).
\ep

Recall that $\bbb_h$ and $\bbb$ denote  the Dirichlet forms associated to $1-T_h$ and $L_\rho$ respectively.
 One has the following
\begin{lem}\label{lem:convDirichform}
Let $\Omega$ be an open set satisfying Assumption \ref{hyp1}. Suppose that 
$(u_h)_{h\in]0,1]}$ is a bounded family of functions in $L^2(\Omega)$ and assume there exists a decomposition 
$u_h=\varphi_h+v_h+r_h$  such that the following assumptions hold true
\begin{itemize}
\item[-] $(\varphi_h)$ converges weakly in $H^1(\Omega)$ towards a limit $\varphi$ when $h\rightarrow 0$.
\item[-] $\Vert v_h\Vert_{L^2}=\ooo(h)$ when $h\rightarrow 0$.
\item[-] $\supp(r_h)\subset \Omega_{c_0 h}$ for some $c_0>0$
\item[-] there exists $C>0$ such that 
$
\eee_h(r_h)\leq C h^2
$
for all $h\in]0,1]$.
\end{itemize}
Then for all $\theta\in H^1(\Omega)$, one has
$$
\lim_{h\rightarrow 0}h^{-2}\bbb_h(u_h,\theta)=\bbb(\varphi,\theta).
$$
\end{lem}

\bp Let us denote $\widetilde \bbb_h=h^{-2}\bbb_h$ and let $\theta\in H^1(\Omega)$. 
We have to prove that
\begin{enumerate}
\item[i)] $\lim_{h\rightarrow 0}\widetilde\bbb_h(r_h,\theta)=0$
\item[ii)]  $\lim_{h\rightarrow 0}\widetilde\bbb_h(\varphi_h,\theta)=\bbb(\varphi,\theta)$
\item[iii)]  $\lim_{h\rightarrow 0}\widetilde\bbb_h(v_h,\theta)=0$.
\end{enumerate}
Let $M\geq 1$ denote a parameter to be fixed later and let $\Omega_{Mh}^c=\Omega\setminus\Omega_{Mh}$ with $\Omega_{Mh}$  defined by 
\eqref{eq:defOmegar}.
Given two subset $A,B$ of $\Omega$ we denote 
$$
\widetilde\bbb_h^{A,B}(u,v)=\frac 1{2 h^{d+2}\Vr_d}\iint_{x\in A,y\in B}\indic_{\vert x-y\vert<h}(u(x)-u(y))\overline{(v(x)-v(y))}d\mu_\rho^2(x,y)
$$
and when $A=B$ we denote $\widetilde\bbb_h^{A,B}(u,v)=\widetilde\bbb_h^{A}(u,v)$. By Cauchy-Schwarz inequality, one has
\be\label{eq:CSDirich}
\widetilde\bbb_h^{A,B}(u,v)\leq h^{-2}\sqrt{\eee_h^A(u)\eee_h^B(v)}.
\ee
Since $(\varphi_h)$ is bounded in $H^1$, it follows from Corollary \ref{cor:locEh} that $\eee_h(\varphi_h)=O(h^2)$. On the other hand,  $\Vert v_h\Vert_{L^2}=O(h)$ implies $\eee_h(v_h)=O(h^2)$  and hence
$\eee_h(u_h)=O(h^2)$. Suppose now that $u_h$, $\theta$ are as above. We claim that 

\be\label{eq:locdirich0}
\widetilde\bbb_h(u_h,\theta)=\widetilde\bbb_h^{\,\Omega^c_{Mh}}(u_h,\theta)+o(1).
\ee
Indeed, one has
\begin{equation*}
\begin{split}
\widetilde\bbb_h(u_h,\theta)=\widetilde\bbb_h^{\,\Omega_{Mh},\Omega_{Mh}}(u_h,\theta)&+\widetilde\bbb_h^{\,\Omega^c_{Mh},\Omega_{Mh}}(u_h,\theta)+\widetilde\bbb_h^{\,\Omega_{Mh},\Omega^c_{Mh}}(u_h,\theta)\\
&+\widetilde\bbb_h^{\,\Omega^c_{Mh},\Omega^c_{Mh}}(u_h,\theta),
\end{split}
\end{equation*}
 and it follows from \eqref{eq:CSDirich} that 
\be\label{eq:locdirich1}
\begin{split}
|\widetilde\bbb_h(u_h,\theta)-&\widetilde\bbb_h^{\,\Omega^c_{Mh},\Omega^c_{Mh}}(u_h,\theta)|\leq 
h^{-2}\sqrt{\eee_h^{\Omega_{Mh}}(u_h)\eee_h^{\Omega_{Mh}}(\theta)}\\
&\phantom{*****}+h^{-2}
\sqrt{\eee_h^{\Omega^c_{Mh}}(u_h)\eee_h^{\Omega_{Mh}}(\theta)}
+|\widetilde\bbb_h^{\,\Omega_{Mh},\Omega^c_{Mh}}(u_h,\theta)|
\\\
&\phantom{*****}\leq 2h^{-2}\sqrt{\eee_h^{}(u_h)\eee_h^{\Omega_{Mh}}(\theta)}+|\widetilde\bbb_h^{\,\Omega_{Mh},\Omega^c_{Mh}}(u_h,\theta)|.
\end{split}
\ee
Since the operator $T_h$ localizes at scale $h$ 
one has 
$$
|\widetilde\bbb_h^{\,\Omega_{Mh},\Omega^c_{Mh}}(u_h,\theta)|=|\widetilde\bbb_h^{\,\Omega_{Mh},\Omega^c_{Mh}\cap\Omega_{(M+1)h}}(u_h,\theta)|
\leq h^{-2}\sqrt{\eee_h(u_h)\eee_h^{\Omega_{(M+1)h}}(\theta)}.
$$
Combining this estimate with \eqref{eq:locdirich1} and using the bound $\eee_h(u_h)=O(h^2)$, we obtain
$$
|\widetilde\bbb_h(u_h,\theta)-\widetilde\bbb_h^{\,\Omega^c_{Mh},\Omega^c_{Mh}}(u_h,\theta)|\leq Ch^{-1}\sqrt{\eee_h^{\Omega_{(M+1)h}}(\theta)}.
$$
By Corollary \ref{cor:locEh}, one knows that $\eee^{\Omega_{(M+1)h}}_h(\theta)=o(h^2)$ which proves \eqref{eq:locdirich0}.

Since $\eee_h(r_h)=O(h^2)$, \eqref{eq:locdirich0} implies 
$\widetilde\bbb_h(r_h,\theta)=\widetilde\bbb_h^{\Omega^c_{Mh}}(r_h,\theta)+o(1)$ and since for $M>c_0$, $r_h=0$ on $\Omega^c_{Mh}$, we get directly i).

Let us now prove ii). Using a partition of unity, we can write $\theta=\sum_{j\in J}\theta_j$ with $\theta_j$ supported in $\omega_j$ for all $j\in J$. Since both side of the equality in ii) are linear with respect to $\theta$ we can assume from now that $\theta$ is supported in a small chart 
$\omega_j$.
Using the change of variable $y=x+hz$, one has 
\begin{equation*}
\begin{split}
\widetilde\bbb_h^{\Omega_{Mh}^c}(\varphi_h,\theta)=
\frac 1{2 h^2 \Vr_d}\int_{\Omega_{Mh}^c}\int_{z\in D_{x,h}}&(\varphi_h(x)-\varphi_h(x+hz))\\
&(\theta(x)-\theta(x+hz))w_h(x,z)dzdx
\end{split}
\end{equation*}
where $D_{x,h}=\{z\in \R^d,\;|z|<1\text{ and }x+hz\in\Omega\}$. Since $M\geq 1$, for any $x\in \Omega_{Mh}^c$, one has $D_h=D=\{|z|<1\}$ and 
for any $t\in[0,1]$ we get $x+thz\in\Omega$. Using this path 
and the argument of Lemma \ref{lem:estimEhUH1}, we can write
\begin{equation*}
\begin{split}]\frac 12,1[
\widetilde \bbb_h^{\Omega_{Mh}^c}(\varphi_h,\theta)=\frac 1{2 \Vr_d}\int_{\Omega_{Mh}^c}&\int_{\vert z\vert<1}
\Big(\int_0^1  z\cdot \nabla\varphi_h( x+shz))ds\Big)\\
&\Big(\int_0^1z\cdot\nabla\theta(x+thz)dt\Big)w_h(x,z)dzdx.
\end{split}
\end{equation*}
 Since $\rho$ is $\ccc^1$, then $w_h(x,z)=\rho(x)+\ooo(h)$ and hence
\begin{equation*}
\begin{split}
\widetilde \bbb_h^{\Omega_{Mh}^c}(\varphi_h,\theta)=\frac 1{2 \Vr_d}\int_{\vert z\vert<1}&\int_{\Omega_{Mh}^c}
\Big(\int_0^1  z\cdot \nabla\varphi_h( x+shz))ds\Big)\\
&\Big(\int_0^1z\cdot\nabla\theta(x+thz)dt\Big)\rho(x)dxdz+O(h).
\end{split}
\end{equation*}
Using the change of variable $\kappa_{s,h,z}:x\mapsto x-shz$, this implies
\begin{equation}\label{eq:bilin00}
\begin{split}
\widetilde \bbb_h^{\Omega_{Mh}^c}(\varphi_h,\theta)=\frac 1{2  \Vr_d}\int_0^1&\int_0^1\int_{\vert z\vert<1}\int_{x\in V_{s,h,z}}
\Big(z\cdot \nabla\varphi_h( x))\Big)\\
&\Big(z\cdot\nabla\theta(x+(t-s)hz)\Big)\rho(y)dxdzdsdt+O(h)
\end{split}
\end{equation}
where $V_{s,h,z}=\kappa_{s,h,z}^{-1}(\Omega_{Mh}^c)$.
We claim that 
\be\label{eq:bilin0}
\int_0^1\int_0^1\int_{\vert z\vert<1}\int_{V_{s,h,z}}\vert \nabla\theta(x+(t-s)hz)-\nabla\theta(x)\vert^2dxdzdtds=o(1).
\ee
Indeed by density of $C^\infty(\Omega)\cap H^1(\Omega) $ in  $ H^1(\Omega)$ we can assume that 
$\theta\in C^\infty(\Omega)$. 
Let us fix $\epsilon>0$. Since $\theta\in H^1(\Omega)$, there exists $r>0$ such  $\int_{\Omega_r}\vert \nabla\theta(x)\vert^2dx\leq \epsilon^2$.
Moreover, since $\nabla\theta$ is uniformly continuous on 
$\overline \Omega_{r}^c$, there exists $h_0>0$ such that for all $h\in]0,h_0]$
$$
\int_0^1\int_0^1\int_{\vert z\vert<1}\int_{\Omega_r^c\cap V_{s,h,z}}\vert \nabla\theta(x+(t-s)hz)-\nabla\theta(x)\vert^2dxdzdtds<\epsilon^2.
$$
Combining these two estimates, we obtain  \eqref{eq:bilin0}. Combined to  \eqref{eq:bilin00}, it implies
\begin{equation*}
\begin{split}
\widetilde \bbb_h^{\Omega_{Mh}^c}(\varphi_h,\theta)&=\frac 1{2  \Vr_d}\int_0^1\int_{\vert z\vert<1}\int_{y\in \Omega}
\Big(z\cdot \nabla\varphi_h( y))\Big)\\
&\phantom{*********}\Big( \indic_{V_{s,h,z}}(y)z\cdot\nabla\theta(y)\Big)\rho(y)dydzds+o(1)
\end{split}
\end{equation*}
Moreover, since $\varphi_h$ is bounded in $H^1$, $\theta\in H^1$ and $\indic_{\Omega\setminus V_{s,h,z}}\rightarrow 0$ pointwise, it follows from
Cauchy-Schwarz inequality and dominated convergence theorem that 
$$
\widetilde \bbb_h^{\Omega_{Mh}^c}(\varphi_h,\theta)=\frac 1{2  \Vr_d}\int_{\vert z\vert<1}\int_{y\in\Omega}
\Big(z\cdot \nabla\varphi_h( y))\Big)\Big(z\cdot\nabla\theta(y)\Big)\rho(y)dydz+o(1)
$$ 
and since $\varphi_h$ converges weakly to $\varphi$ in $H^1$, we obtain
\begin{equation*}
\begin{split}
\lim_{h\rightarrow 0}\widetilde\bbb_h^{\Omega_{Mh}^c}(\varphi_h,\theta)&=\frac 1{2 \Vr_d}\int_\Omega\int_{\vert z\vert<1}(z\cdot\nabla\varphi(x))(z\cdot\nabla\theta(x))dz\rho(x)dx\\
&=\frac 1{2 \Vr_d}\sum_{i,j=1}^d\int_\Omega\int_{\vert z\vert<1}z_i\partial_i\varphi(x)z_j\partial_j\theta(x)dz\rho(x)dx.
\end{split}
\end{equation*}
For parity reason the terms associated to $i\neq j$ vanish and using \eqref{eq:locdirich0}, we get
$$
\lim_{h\rightarrow 0}\widetilde\bbb_h(\varphi_h,\theta)=\lim_{h\rightarrow 0}\widetilde\bbb_h^{\Omega_{Mh}^c}(\varphi_h,\theta)=
\sum_{i=1}^da_j\int_\Omega\partial_i\varphi(x)\partial_i\theta(x)\rho(x)dx
$$
with 
$$a_j=\frac 1{2 \Vr_d}\int_{\vert z\vert<1}z_i^2dz=\frac 1{2(d+2)}.$$
This proves ii).

It remains to prove iii). As before we can work with the functional $\widetilde\bbb_h^{\Omega_{Mh}^c}$ instead of $\widetilde\bbb_h$.
 One has 
\begin{equation*}
\begin{split}
\widetilde\bbb_h^{\Omega_{Mh}^c}(v_h,\theta)=\frac 1{2 h^2\Vr_d}\int_\Omega\int_{\vert z\vert<1}&(v_h(x)-v_h(x+hz))\\
&(\theta(x)-\theta(x+hz) w_h(x,z)dzdx
\end{split}
\end{equation*}
Splitting the difference $(v_h(x)-v_h(x+hz))$ in two different integrals and making the change of variable $x\mapsto x-hz$ in the term corresponding to 
$v_h(x+hz)$ we get
$$
\widetilde\bbb_h^{\Omega_{Mh}^c}(v_h,\theta)= \widetilde \bbb_h^+(v_h,\theta)+\widetilde \bbb_h^-(v_h,\theta)
$$
with 
$$\bbb_h^\pm(v_h,\theta)=\frac 1{2 h^2 \Vr_d}\int_{\Omega_{Mh}^c}\int_{\vert z\vert<1}v_h(x)(\theta(x)-\theta(x\pm hz))
 w_h(x,z)dzdx,
$$
where the integration domain in the variable $z$ is the unit disc for the same reason as before.
We show how to estimate $\widetilde\bbb_h^+$, the case of $\widetilde\bbb_h^-$ is similar.
The same computation as above shows that 
\begin{equation*}
\begin{split}
\widetilde\bbb_h^+(v_h,\theta)=\frac 1{2 h \Vr_d}&\int_{\vert z\vert<1}\int_{\Omega_{Mh}^c}v_h(x)\\
&\Big(\int_0^1z\cdot 
\nabla\theta(x+thz)dt\Big) \rho(x)dxdz +O(h)
\end{split}
\end{equation*}
where we used again  $w_h(x,z)=\rho(x)+O(h)$. Since $\theta\in H^1$,
$\Vert v_h\Vert_{L^2}=O(h)$ and $\indic_{\Omega\setminus \Omega_{Mh}^c}\rightarrow 0$ pointwise,  we get as in the proof of ii) that
$$
\widetilde\bbb_h^+(v_h,\theta)=\frac 1{2 h\Vr_d}\int_{\Omega}\int_{\vert z\vert<1}v_h(x)\Big(z\cdot 
\nabla\theta(x)\Big) \rho(x)dzdx +o(1),
$$
and since $\int_{\vert z\vert<1}(z\cdot 
\nabla\theta(x))dz=0$, we obtain $\widetilde\bbb_h^+(v_h,\theta)=o(1)$ which proves iii).
\ep

\subsection{Case of smooth densities}
In this section we prove  Theorem \ref{th:spec-reg}.
We follow the proof of Theorem 1.2 in 
\cite{DiLeMi12}. 
Let $\vert \triangle_{h}\vert $ be the rescaled (non negative)
Laplacien associated to the Markov kernel $T_{h}$
\be\label{5.-1}
\vert \triangle_{h} \vert ={1-T_{h}\over h^2}.
\ee
Let $R>0$ be fixed. If
$\nu_h\in[0,R]$ and $u_h\in L^2(M)$ satisfy $| \triangle_h|
u_h=\nu_h u_h$ and $\| u_h\|_{L^2}=1$, then thanks to Lemma
\ref{lem:decompHF}, $u_h$ can be decomposed as $u_h=\varphi_h+v_h+r_h$ with
$\|v_h\|_{L^2}=O(h)$, $\varphi_h$ bounded in $ H^1(\Omega)$ and $r_h$ supported in $\Gamma_h\subset\Omega_{c_0h}$ for some $c_0>0$.
Moreover, we claim that $\eee_h(r_h)=O(h^2)$. Indeed, since $r_h=u_h-\varphi_h-v_h$ and $\eee_h(u_h)=h^2\nu_h$ with $\nu_h$ bounded,  it suffices to show that 
$\eee_h(v_h)$ and $\eee_h(\varphi_h)$ are $O(h^2)$. The bound on $\eee_h(v_h)$ follows directly from the fact that $\Vert r_h\Vert_{L^2}=\ooo(h)$ and that 
$1-T_h$ is bounded on $L^2$. The bound on $\eee_h(\varphi_h)$ is obtained from the fact that $\varphi_h$ is bounded in $H^1$ and Corollary \ref{cor:locEh}.
Consequently, 
(extracting a subsequence if necessary)  we can assume that
$(\varphi_h)$ weakly converges in $ H^1(\Omega)$ to a limit $\varphi$ and that $(\nu_h)$
converges to a limit $\nu$. Hence $(u_{h})$ converge strongly in $L^2$ to $\varphi$, and 
it now follows from Lemma \ref{lem:convDirichform} that
for any $\theta\in C^\infty(M)$,
\be\label{T.40}
\begin{split}
\nu\<\varphi,\theta\>=\lim_{h\rightarrow 0}\nu_h\<u_h,\theta\>=\lim_{h\rightarrow 0} h^{-2}\bbb_h(u_h,\theta)= \bbb(\varphi,\theta).
\end{split}
\ee
Since $\theta$ is arbitrary in $H^1$ this shows that $\varphi\in D(L_\rho)$ and that $(L_\rho-\nu)\varphi =0$. 
Hence $\nu$ is an eigenvalue of $L_\rho$.
Moreover, the dimension of an orthonormal basis is preserved by strong limit. So the above argument   proves that
for any $\epsilon>0$ small, there exists $h_\epsilon>0$ such that for $h\in]0,h_\epsilon]$, one has
\begin{equation}
\label{T41}
\sigma(\vert\Delta_h\vert)\cap[0,R]\subset\cup_j[\nu_j-\epsilon,\nu_j+\epsilon]
\end{equation}
and
\begin{equation}
\sharp \sigma(|\Delta_h|)\cap[\nu_j-\epsilon,\nu_j+\epsilon]\leq m_j.
\label{T42}
\end{equation}
In order to show that one has equality in  \eqref{T42} for $\epsilon$ small enough, observe that 
for any $\psi\in H^1(\Omega)$ independent of $h$, one has 
$$
\lim_{h\rightarrow 0}h^{-2}\eee_h(\psi)=\bbb(\psi,\psi)
$$
thanks to Lemma \ref{lem:convDirichform}.
In particular, if $\psi\in D(L_\rho)$ satisfies $L_\rho\psi=\nu\psi$ for some $\nu>0$, then 
 $\lim_{h\rightarrow 0} h^{-2}\eee_h(\psi)=\nu \Vert \Psi\Vert^2$.
Hence, we can mimic the proof 
of Theorem 2 iii) in \cite{DiLeMi12} to get the result. The proof of Theorem \ref{th:spec-reg} is complete.

\subsection{Case of measurable densities.}
In this section we assume that $\rho$ is a measurable function satisfying \eqref{eq:bornerho} and we prove Theorem \ref{th:spec-rough}. We first apply Theorem \ref{th:spec-reg} with $\rho_0=1$.
It follows that $1$ is a simple eigenvalue of $T_{h,\rho_0}$. Moreover, denoting $(\mu_{k,\rho_0}(h))_{k\in\N}$ the decaying sequence of positive eigenvalues of $T_{h,\rho_0}$, one has $1=\mu_{0,\rho_0}>\mu_{1,\rho_0}(h)$ and $\mu_{1,\rho_0}(h)=h^2\nu_1+o(h^2)$ where we recall that 
$\nu_1>0$ is the lowest positive eigenvalue of the Neumann Laplacian on $\Omega$.
Moreover, one has $\ker(T_{h,\rho_0}-1)=\Span(1)$
Combined to the spectral theorem, this implies that  for all $u\in\Span(1)^\bot$, we have
\be\label{eq:pfspecr1}
\<(1-T_{h,\rho_0})u,u\>_{L^2(\rho_0)}\geq Ch^2\Vert u\Vert^2_{L^2(\rho_0)}.
\ee
On the other hand, from \eqref{eq:defDirich1} one has 
$$
\<(1-T_{h,\rho_0})u,u\>_{L^2(\rho_0)}=\frac 1{2 h^{d}\Vr_d}\int_{\Omega\times \Omega}\indic_{\vert x-y\vert<h}(f(x)-f(y))^2 d\mu_{\rho_0}^2(x,y),
$$
and since $m\leq \rho\leq M$, then 
$$
\<(1-T_{h,\rho_0})u,u\>_{L^2(\rho_0)}\leq \frac 1 m\<(1-T_{h,\rho})u,u\>_{L^2(\rho)}
$$
and $\Vert u\Vert^2_{L^2(\rho_0)}\geq \frac 1M\Vert u\Vert^2_{L^2(\rho)}$. Combined with \eqref{eq:pfspecr1}, this implies that there exists a
new positive constant $C$ such that  for all $u\in\Span(1)^\bot$, we have
$$
\<(1-T_{h,\rho})u,u\>_{L^2(\rho)}\geq Ch^2\Vert u\Vert^2_{L^2(\rho)}.
$$ 
This proves i) and the lower bound on $g(h)$. The upper bound is proved in the same way, using the equivalence of Dirichlet forms.
\subsection{Total variation estimates.} 
This section is devoted to the proof of Theorem \ref{th:TVestim}.
Thanks to \eqref{eq:TVdist1}, we have
$$
\sup_{x\in\Omega}\Vert t_{h,\rho}^n(x,dy)-\mu_\rho\Vert_{TV}=\frac 12\Vert T_{h,\rho}^n-\Pi_0\Vert_{L^\infty\rightarrow L^\infty}
$$
where $\Pi_0$ denotes the orthogonal projection on $\Span(1)$ in $L^2(\rho)$. Throughout this section, we drop the dependance with respect to $\rho$ in the notations. For any $p\in\N$, one has $T_h^p=A_p+B_p$ with $A_1=m_h$, $B_1=K_h$ and for any $p\geq 1$
$A_{p+1}=m_h A_p$, $B_{p+1}=m_hB_p+K_h T_h^p$.
Since $\Vert m_h\Vert_{L^\infty\rightarrow L^\infty}\leq 1-Ch^\gamma$ and $\Vert K_h\Vert_{L^2\rightarrow L^\infty}\leq Ch^{-\frac d2}$, 
it follows from (2.49) and (2.50) in \cite{DiLeMi11} that for any $p\in\N$
\be\label{eq:pfTVestim1}
\begin{split}
\Vert A_p\Vert_{L^\infty\rightarrow L^\infty}&\leq(1-Ch^\gamma)^p\\
\Vert B_p\Vert_{L^2\rightarrow L^\infty}&\leq C h^{-\gamma-\frac d 2}.
\end{split}
\ee
Suppose now that $p,n\in\N$. Since $T_h\Pi_0=\Pi_0$ we get
\begin{equation*}
\begin{split}
\Vert T_h^{p+n+1}-\Pi_0\Vert_{L^\infty\rightarrow L^\infty}\leq 
\Vert A_p\Vert_{L^\infty\rightarrow L^\infty}&\Vert T_h^{n+1}-\Pi_0\Vert_{L^\infty\rightarrow L^\infty}\\
&+\Vert B_p(T_h^{n+1}-\Pi_0)\Vert_{L^\infty\rightarrow L^\infty}
\end{split}
\end{equation*}
Taking $p=\lfloor M n h^{2-\gamma}\rfloor$ with $M>0$ to be chosen large enough (here we denote $\lfloor n\rfloor$ the integer part of $n\in\N$), we deduce from \eqref{eq:pfTVestim1} that 
$$
\Vert A_p\Vert_{L^\infty\rightarrow L^\infty}\leq e^{-nMCh^2}
$$
where $C$ is a positive constant independent of $h$ and $M$.
Since $T_h$ is markovian, $T_h$ and $\Pi_0$ are bounded by $1$ on $L^\infty$ and conseqently
\be\label{eq:pfTVestim2}
\Vert T_h^{p+n+1}-\Pi_0\Vert_{L^\infty\rightarrow L^\infty}\leq C e^{-nMCh^2}
+\Vert B_p(T_h^{n+1}-\Pi_0)\Vert_{L^\infty\rightarrow L^\infty}.
\ee
We shall now estimate the second term in the above right hand side. One has 
$$
\Vert B_p(T_h^{n+1}-\Pi_0)\Vert_{L^\infty\rightarrow L^\infty}\leq 
\Vert B_p\Vert_{L^2\rightarrow L^\infty}\Vert T_h^n-\Pi_0\Vert_{L^2\rightarrow L^2}
\Vert T_h\Vert_{L^\infty\rightarrow L^2}
$$
and from Proposition \ref{prop:spec-ess-Th} and Theorem \ref{th:spec-rough}, we know that 
$\sigma(T_h)\setminus\{1\}\subset [-1+Ch^\gamma,1-g(h)]$ with $h^2/C\leq g(h)\leq C h^2$ and $\gamma<2$. Hence it follows from 
the spectral theorem, that for $h$ small enough
$$
\Vert T_h^n-\Pi_0\Vert_{L^2\rightarrow L^2}\leq (1-g(h))^n.
$$
Combined with \eqref{eq:pfTVestim1} and the estimate 
$\Vert T_h\Vert_{L^\infty\rightarrow L^2}\leq \Vert T_h\Vert_{L^\infty\rightarrow L^\infty}=1$, it follows that
$$
\Vert B_p(T_h^{n+1}-\Pi_0)\Vert_{L^\infty\rightarrow L^\infty}
\leq 
C h^{-\gamma-\frac d 2} (1-g(h))^n\leq C h^{-\gamma-\frac d 2} e^{-ng(h)}.
$$
Together with \eqref{eq:pfTVestim2}, this implies
$$
\Vert T_h^{p+n+1}-\Pi_0\Vert_{L^\infty\rightarrow L^\infty}\leq C e^{-nMCh^2}+C h^{-\gamma-\frac d 2} e^{-ng(h)}.
$$
Since $h^2/C\leq g(h)\leq C h^2$, it follows that for $M>0$ large enough one has
$$
\Vert T_h^{p+n+1}-\Pi_0\Vert_{L^\infty\rightarrow L^\infty}\leq C h^{-\gamma-\frac d 2} e^{-ng(h)}.
$$
Taking advantage of $p=\lfloor Mn h^{2-\gamma}\rfloor$, this can be written
$$
\Vert T_h^{n}-\Pi_0\Vert_{L^\infty\rightarrow L^\infty}\leq C h^{-\gamma-\frac d 2} e^{-ng(h)(1+O(h^{2-\gamma}))}
$$
which proves \eqref{1.7}.

\section{Appendix}\label{A2}
Let $\Pi^d=(\R/2\Z)^d$, $d=1+d'+d''$.
 For $f\in L^2(\Pi^d)$ and for any $k=(k_1,k',k'')\in\Z\times\Z^{d'}\times \Z^{d''}$, we denote by
 \be\label{eq:defFouriersc}
\begin{split}
\hat f(k)=\fff f(k):=
\frac 1{2^{d/2}}\int_{\Pi^d}&e^{-i\pi\<x,k\> }f(x) dx.
\end{split}
\ee
the Fourier coefficients of the function $f$.
The map $\fff$ is  an isometry from $L^2(\Pi^{d})$ onto $\ell^2(\Z^{d})$ and we denote by $\bar \fff$ its adjoint:
\be
\bar \fff(a)=\frac 1{2^{d/2}}\sum_{k\in\Z^d}a_k e^{i\pi\<x,k\>}
\ee
for any $a=(a_k)_{k\in\Z^d}$.
Let also $\lambda_1,\lambda',\lambda''>0$ be some parameters and denote $\lambda=(\lambda_1,\lambda',\lambda'')$. We recall that 
for any $\xi=(\xi_1,\xi',\xi'')\in \R^{1+d'+d''}$ we denote $\lambda\cdot \xi=(\lambda_1 \xi_1,\lambda'\xi',\lambda''\xi'')$.
For any $s\in\R$, we define  the $\lambda$-Sobolev space as the space of functions $\phi$ such that $\Vert \phi\Vert_{H^s_{\lambda}}<\infty$ where
\be\label{eq:definormSobsc}
\Vert \phi\Vert_{H^s_{\lambda}}=\Vert(\<\lambda\cdot k\>^s\fff f(k))_k\Vert_{\ell^2(\Z^d)}.
\ee
We define similarly the partial Fourier coefficients $\fff_{x'x''}:L^2(\Pi^{1+d'+d''})\rightarrow \ell^2(\Z^{d'+d''},L^2(\Pi))$, 
$\fff_{x_1}: L^2(\Pi^{1+d'+d''})\rightarrow \ell^2(\Z,L^2(\Pi^{d'+d''}))$ and their adjoint $\bar \fff_{x',x''}$, $\overline \fff_{x_1}$.
Consider the hypersurface $\Sigma_a=\{x_1=a\}\times\Pi^{d-1}\subset\Pi^d$. We define the trace operator 
$\gamma^\Pi_a: H_\lambda^1(\Pi^d)\rightarrow H^{1/2}_\lambda(\Sigma_a)$ by 
\be\label{eq:defgam0}
\begin{split}
\gamma_a^\Pi \phi(x',x'')&=\frac 1{\sqrt{2 }}\bar \fff_{x',x''}\Big(
\sum_{k_1\in\Z}e^{i\pi k_1a}\fff \phi(k_1,k',k'')\Big)(x',x'')\\
&=\frac 1{\sqrt{2 }}\sum_{k_1\in\Z}e^{i\pi k_1a}\fff_{x_1}\phi(k_1,x',x'').
\end{split}
\ee
\begin{lem}\label{lem:trace-sob} Let $s>\frac 12$. There exists $C>0$ such that for any $\lambda_1,\lambda',\lambda''>0$ and any $\phi\in H^s_\lambda(\Pi^d)$ such that $\fff_{x_1}\phi(0,x',x'')=0$, one has
\be\label{eq:estimatrace}
\Vert\gamma_a^\Pi\phi\Vert_{H^{s-\frac 12}_{\lambda',\lambda''}(\Sigma_a)}\leq C\lambda_1^{-\frac 12}\Vert\phi\Vert_{H^s_\lambda(\Pi^d)}.
\ee
\end{lem}
\bp
We may assume without loss of generality that $a=0$.
 By a density argument, it is sufficient 
to prove \eqref{eq:estimatrace} for $\phi\in C^\infty(\Pi^d)$ such that $\fff_{x_1}\phi(0,x',x'')=0$. 
For such functions, the sum in  \eqref{eq:defgam0} is over $k_1\in \Z^*$ and it follows from Cauchy-Schwarz inequality that
\be\label{eq:trace-sob1}
\begin{split}
|\fff_{x',x''}&(\gamma_a \phi)(k',k'')|=\frac 1{\sqrt{2 }}|\sum_{k_1\in\Z^*}e^{i\pi k_1 a}\fff \phi(k)|\\
&\leq \frac 1{\sqrt{2 }}\Big( \sum_{k_1\in\Z^*}\<\lambda\cdot k\>^{-2s} \Big)^{\frac 12}
\Big( \sum_{k_1\in\Z^*}\<\lambda\cdot k\>^{2s}|\fff\phi(k)|^2 \Big)^{\frac 12}.
\end{split}
\ee
We claim that there exists a constant $C>0$ such that for any $(k',k'')\in \Z^{d'+d''}$ and any $\lambda',\lambda''>0$, one has
\be\label{eq:trace-sob2}
\sum_{k_1\in\Z^*}\<\lambda\cdot k\>^{-2s}\leq C\lambda_1^{-1}\<(\lambda' k',\lambda'' k'')\>^{1-2s}.
\ee
Indeed, since the function $m:t\mapsto (1+|\lambda_1 t|^2+|\lambda' k'|^2+|\lambda'' k''|^2)^{-s}$ is decreasing and integrable on $\R$, one has 
\be\label{eq:trace-sob3}
\begin{split}
\sum_{k_1\in\Z^*}\<\lambda\cdot k\>^{-2s}&=\sum_{k_1\in\Z^*} m(k_1)\leq \int_\R m(t)dt.
\end{split}
\ee
Using the change of variable $t\mapsto \frac {\<(\lambda' k',\lambda'' k'')\>}{\lambda_1}t$ one gets
$\int_\R m(t)dt=C_1\lambda_1^{-1}\<(\lambda' k',\lambda'' k'')\>^{1-2s}$ for some universal constant $C_1$. Combined with \eqref{eq:trace-sob3}, this proves  \eqref{eq:trace-sob2}.

Now, using  \eqref{eq:trace-sob2} and  \eqref{eq:trace-sob1}, we get
$$
|\fff_{x',x''}(\gamma_a \phi)(k',k'')|^2\leq C\lambda_1^{-1}\<(\lambda' k',\lambda'' k'')\>^{1-2s} \sum_{k_1\in\Z^*}\<\lambda\cdot k\>^{2s}|\fff\phi(k)|^2
$$
and hence
\begin{equation*}
\begin{split}
\Vert \gamma_a^\Pi\phi\Vert_{H^{s-\frac 12}_{\lambda',\lambda''}}^2&= \sum_{k',k''}\<(\lambda' k',\lambda'' k'')\>^{2s-1}|\fff_{x',x''}(\gamma_a \phi)(k',k'')|^2\\
&\leq C\lambda_1^{-1} \sum_{k',k''}\sum_{k_1\neq 0}\<\lambda\cdot k\>^{2s}|\fff\phi(k)|^2=C\lambda_1^{-1}\Vert \phi\Vert_{H^s_\lambda(\Pi^d)}^2
\end{split}
\end{equation*}
which proves the result.
\ep

Given $0<a<b<2$, the restriction operator defined by \eqref{eq:defrestrict} acts on $H^1$ functions 
$
R_{]a,b[}:H^1(\Pi^d)\rightarrow H^1(]a,b[\times\Pi^{d-1})
$ and one defines
the trace operator 
\be\label{eq:deftrace}
\gamma_a:H^1(]a,b[\times\Pi^{d-1})\rightarrow H^{\frac 12}(\Sigma_a)
\ee
by $\gamma_af=\gamma_a^\Pi\tilde f$ for any $\tilde f\in H^1(\Pi^d)$ such that $R_{]a,b[}\tilde f=f$.
Throughout we write $\gamma_af=f_{|x_1=a}$

Suppose now that $a<b<c$ are some fixed real numbers and let
 $A_0,A_1,A_2\subset\R\times\Pi^{d-1}$ be defined by $A_0=]a,b[\times \Pi^{d-1}$,
$A_1=]b,c[\times \Pi^{d-1}$, 
$A_2= ]a,c[\times \Pi^{d-1}$. 
\begin{lem}\label{lem:recol}
Let $(\phi_j)_{j=0,1,2}\in H_\lambda^1 (A_{j})$ and $r_2\in L^2(A_{2})$ be some functions depending one some parameters 
$\lambda=(\lambda_1,\lambda',\lambda'')\in]0,+\infty[^3$ and $h>0$. Let 
$f\in L^2(A_{2})$ 
given by $f=\mathds{1}_{A_{0}}\phi_0+ \mathds{1}_{A_{1}}\phi_1$ and assume that $f=\phi_2+r_2$ with
$$\Vert \phi_j\Vert_{H^1_\lambda(A_j)}\leq 1 \text{ and } 
\Vert r_2\Vert_{L^2(A_{2})}\leq h$$
for all $j=0,1,2$.
 Then there exists $h_1>0$  and $\Upsilon>0$ such that for  $0<\lambda_1 h<h_1$,  there exists 
 $\psi\in H^1_\lambda(A_1)$ supported in $b \leq x_1<b+h\lambda_1 $ and such that
$\psi_{\lvert x_1=b}=( \phi_0)_{\lvert x_1=b}-(\phi_1)_{\lvert x_1=b}$ and  
\be\label{deux}
 \Vert\psi\Vert_{H^1_\lambda(A_1)}\leq \Upsilon\text{ and }\lvert \lvert \psi\rvert \rvert_{L^2(A_{1})} \leq  \Upsilon  h
\ee
\end{lem}

\bp Throughout $C$ denotes a positive constant independent of $h$ and $\lambda$ that may change from line to line. First observe that the statement of the lemma is invariant by translation and dilation in the variable $x_1$. Hence we can assume without loss of generality that $a=-1,b=0$ and $c=1$. Throughout the proof, we  denote  $\Sigma=\{x_1=0\}\times\Pi^{d-1}\subset ]a,c[\times\Pi^{d-1}$ and 
we let $\sigma:]a,c[\times\Pi^{d-1}\rightarrow ]a,c[\times\Pi^{d-1}$ denote the symmetry with respect to $\Sigma$. We define 
$g_0=\mathds{1}_{A_{0}}\phi_0+ \mathds{1}_{A_{1}}\phi_0\circ\sigma$ and 
$g_1=\mathds{1}_{A_{0}}\phi_1\circ\sigma+ \mathds{1}_{A_{1}}\phi_1$. We denote 
\be\label{eq:deftheta}
\theta= ( \phi_0)_{|\Sigma}-( \phi_1)_{|\Sigma}:=\gamma_0(g_0-g_1)
\ee
with $\gamma_0$ defined by \eqref{eq:deftrace}.
We claim that there exists $C_0>0$ independent of the $\phi_i$ such that
\begin{equation}\label{eq:estimthetaL2}
\Vert \theta \Vert_{L^2(\Sigma)} \leq C_0 \sqrt{\frac h{\lambda_1}}.
\end{equation}
In order to prove this estimate, let $\varepsilon>0$ a constant to be fixed later and let 
\be\label{eq:defIeps}
I_\epsilon=I_\epsilon(x'):=\int_{-\epsilon}^0  f(x_1,x')dx_1- \int_0^\epsilon  f(x_1,x')dx_1
\ee
wich is well defined for $\vert \epsilon\vert<h_1:=\min(b-a,c-b)$.
By Taylor expansion,  one has
$\phi_i(x)=\phi_i(0,x')+ \int_0^{x_1} \partial_1 \phi_i(t,x')dt$
and hence
\begin{equation}\label{eq:expandIeps}
\begin{split}
I_\epsilon(x')&=\int_{-\epsilon}^0 \phi_0(x_1,x')dx_1 -\int_0^\epsilon  \phi_1(x_1,x')dx_1\\
&=\epsilon \theta(x')+\int_{-\epsilon}^0 \int_0^{x_1}\partial_1\phi_0 (t,x')dtdx_1\\
& \phantom{********}-\int_0^\epsilon \int_0^{x_1} \partial_1 \phi_1(t,x')dt dx_1.
\end{split}
\end{equation}
Moreover, one has
\begin{equation}\label{eq:estimIeps1}
\begin{split}
\Big\Vert \int_{-\epsilon}^0 \int_0^{x_1}\partial_1\phi_0 (t,x')&dtdx_1\Big\Vert _{L^2(\Pi^{d-1})}\\
&\leq \Big\Vert \int_{-\epsilon}^0 \sqrt {\vert x_1\vert} \lambda_1^{-1}\Vert \lambda_1\partial_1\phi_0\Vert_{L^2(]a,b[)}dx_1\Big\Vert _{L^2(\Pi^{d-1})}\\
& \leq \lambda_1^{-1}\int_{-\epsilon}^0\sqrt{\vert x_1\vert}\Vert \lambda_1\partial_1 \phi_0\Vert _{L^2(A_0)}dx_1\\
&\leq \frac{2}{3}\epsilon^{\frac{3}{2}}\lambda_1^{-1} \Vert  \phi_0\Vert _{H^1_\lambda(A_0)}\leq C\epsilon^{\frac{3}{2}}\lambda_1^{-1},
\end{split}
\end{equation}
and of course an estimate similar  to \eqref{eq:estimIeps1}  holds true for $\phi_1$.
On the other hand,  since $f=\phi_2+r_2$, then
\begin{equation*}
\begin{split}
I_\epsilon(x')&=\int_{-\epsilon}^0r_2(x_1,x')dx_1-\int^{\epsilon}_0r_2(x_1,x')dx_1+\int_{-\epsilon}^0\phi_2(0,x')dx'
-\int_0^\epsilon \phi_2(0,x')dx'\\
&\phantom{*****}+\int_{-\epsilon}^0 \int_0^{x_1}\partial_1\phi_2 (t,x')dtdx_1 -\int_0^\epsilon \int_0^{x_1} \partial_1 \phi_2(t,x')dt dx_1\\
&=\int_{-\epsilon}^0r_2(x_1,x')dx_1-\int^{\epsilon}_0r_2(x_1,x')dx_1\\
&\phantom{*****}+\int_{-\epsilon}^0 \int_0^{x_1}\partial_1\phi_2 (t,x')dtdx_1 -\int_0^\epsilon \int_0^{x_1} \partial_1 \phi_2(t,x')dt dx_1.
\end{split}
\end{equation*}
The two last terms of the above identity are estimated as above. It follows that 
\begin{equation*}
\Vert I_\epsilon(x')\Vert_{L^2(\Pi^{d-1})}\leq  \int_{-\epsilon}^\epsilon\Big\Vert r_2(x_1,x')\Big \Vert_{L^2(\Pi^{d-1})}dx_1
+C\epsilon^{\frac{3}{2}}\lambda_1^{-1}.
\end{equation*}
Using Cauchy-Schwarz and the assumption on $r_2$, we get 
$$
\Vert I_\epsilon(x')\Vert_{L^2(\Pi^{d-1})}\leq C(h\sqrt \epsilon+\epsilon^{\frac{3}{2}}\lambda_1^{-1}).
$$
Combining this estimate, \eqref{eq:expandIeps} and \eqref{eq:estimIeps1}, we get
$$
\Vert \theta\Vert_{L^2(\Sigma)}\leq C(\frac{\sqrt \epsilon}{\lambda_1}+\frac h{\sqrt\epsilon}).
$$
Minimizing the right hand side by taking  $\epsilon =h\lambda_1$we get $\lvert \lvert \theta\lvert \lvert _{L^2(\Sigma)}=\ooo(\sqrt{h/\lambda_1})$
which proves \eqref{eq:estimthetaL2}. 

Next we want to estimate half derivatives of $\theta$.
We recall that $\theta$ is defined by \eqref{eq:deftheta} and we decompose $g_0-g_1=\delta+\bar g$ with 
\be\label{eq:defgbar}
\bar g(x)=\int_{-1}^1(g_0-g_1)(t,x',x'')dt.
\ee
The function $\bar g$ is independant of $x_1$, hence it is defined as a function on $A_2$ and as a function on $\Sigma$
\begin{slem}\label{slem-trace}
One has $\bar g\in H^1_{\lambda}(A_2)$ and $\bar g\in H^{\frac 12}_{\lambda',\lambda''}(\Sigma)$. Moreover
\be\label{eq:estimgmoyen1}
\Vert \bar g\Vert_{ H^1_{\lambda}(A_2)}\leq \Vert \phi_0\Vert_{ H^1_{\lambda}(A_0)}+\Vert \phi_1\Vert_{ H^1_{\lambda}(A_1)}
\ee
and there exists a constant $C>0$ such that 
\be\label{eq:estimgmoyen2}
\Vert \bar g\Vert_{ H^{\frac 12}_{\lambda',\lambda''}(\Sigma)}\leq \frac C {\sqrt{\lambda_1}}
\ee
\end{slem}
\bp
By Cauchy-Schwarz inequality, one has
$\Vert \bar g\Vert_{L^2(A_2)}\leq C(\Vert \phi_0\Vert_{L^2}+\Vert \phi_1\Vert_{L^2})$,
$\Vert \lambda'\partial_{x'}\bar g\Vert_{L^2(A_2)}\leq C(\Vert \lambda'\partial_{x'}\phi_0\Vert_{L^2}+\Vert \lambda'\partial_{x'}\phi_1\Vert_{L^2})$
 and a similar estimate for derivative in the variable $x''$. Moreover, 
 $\lambda_1\partial_{x_1}\bar g=0$. Hence we have $\bar g\in H^1_{\lambda}(A_2)$ and \eqref{eq:estimgmoyen1} holds true.
By a classical trace theorem, it follows that $\bar g \in H^{\frac 12}_{\lambda',\lambda''}(\Sigma)$ and it remains to prove \eqref{eq:estimgmoyen2}.
 We can assume $\lambda_1\geq 1$.
One has 
$$
\Vert \bar g\Vert_{ H^{\frac 12}_{\lambda',\lambda''}(\Sigma)}^2=\sum_{\tilde k\in\Z^{d'+d''}}\<\tilde\lambda\cdot\tilde k\>\vert\fff_{x',x''}\bar g(\tilde k)\vert^2
$$
where we denote $\tilde\lambda=(\lambda',\lambda'')$ and $\tilde k=(k',k'')$. Splitting the sum in two parts we get
$
\Vert \bar g\Vert_{ H^{\frac 12}_{\lambda',\lambda''}(\Sigma)}^2=S_\leq(\lambda_1)+S_>(\lambda_1)
$
where 
$$
S_>(\lambda_1)=\sum_{\<\tilde\lambda\cdot\tilde k\>> \lambda_1}\<\tilde\lambda\cdot\tilde k\>\vert\fff_{x',x''}\bar g(\tilde k)\vert^2.
$$
One has 
\be\label{eq:estimS>}
S_>(\lambda_1)\leq \frac 1{\lambda_1}\sum_{\<\tilde\lambda\cdot\tilde k\>> \lambda_1}\<\tilde\lambda\cdot\tilde k\>^2\vert\fff_{x',x''}\bar g(\tilde k)\vert^2
\leq\frac 1{\lambda_1}\Vert \bar g\Vert_{ H^1_{\lambda}(A_2)}^2\leq \frac C{\lambda_1}
\ee
thanks to \eqref{eq:estimgmoyen1}.
In order to  estimate the low frequencies, we observe that 
\be\label{eq:estimS<}
S_\leq(\lambda_1)\leq \lambda_1\sum_{\<\tilde\lambda\cdot\tilde k\>\leq \lambda_1}\vert\fff_{x',x''}\bar g(\tilde k)\vert^2\leq \lambda_1\Vert \bar g\Vert_{L^2(\Sigma)}^2.
\ee
We claim that 
\be\label{eq:estimbargL2}
\Vert  \bar g\Vert_{L^2(\Sigma)}^2\leq (h^2+ \frac 1 {\lambda_1^2})
\ee 
Indeed, by Cauchy-Schwarz inequality and thanks to the symetric form of $g_0$ and $g_1$, one has 
\begin{equation*}
\begin{split}
\Vert  \bar g\Vert_{L^2(\Sigma)}&\leq C\Vert\phi_0-\phi_1\circ\sigma\Vert_{L^2(A_0)}=C\Vert f-f\circ\sigma\Vert_{L^2(A_2)}\\
&\leq C\Vert \phi_2-\phi_2\circ\sigma\Vert_{L^2(A_2)}+C\Vert r_h\Vert_{L^2(A_2)}\\
\end{split}
\end{equation*}
Moreover, since $\phi_2(x)-\phi_2\circ\sigma(x)=\int_{-x_1}^{x_1}\partial_1\phi_2(t,x',x'')dt$, we get
\be\label{eq:estimbargL21}
\begin{split}
\Vert  \bar g\Vert_{L^2(\Sigma)}&\leq \frac C{\lambda_1}\Vert\phi_2\Vert_{H^1_\lambda}+C\Vert r_h\Vert_{L^2(A_2)}\leq  \frac C{\lambda_1}+Ch
\end{split}
\ee
which proves \eqref{eq:estimbargL2}. Now combining \eqref{eq:estimS<} and \eqref{eq:estimbargL2} we get 
$ 
S_\leq(\lambda_1)\leq \frac C{\lambda_1}+Ch
$
which combined with \eqref{eq:estimS>} proves the result since $h\lambda_1$ is bounded.
\ep
We are now in position to estimate $\theta$ in $H^{\frac 12}$.
One has $\theta =\bar g+\gamma_0(\delta)$ with $\delta=g_0-g_1-\bar g$ and from Sub-lemma \ref{slem-trace} we know that 
\be\label{eq:estimethetaH012}
\Vert \bar g\Vert_{H^{1/2}_{\lambda',\lambda''}(\Sigma)}\leq C/\sqrt{\lambda_1}.
\ee
Moreover,
by construction, one has 
$\Vert g_j\Vert_{H^{1}_{\lambda}(A_2)}\leq C$ for $j=0,1$ and by Sub-lemma \ref{slem-trace} one has 
also $\Vert \bar g\Vert_{H^{1}_{\lambda}(A_2)}\leq C$. Hence 
$
\Vert \delta\Vert_{H^{1}_{\lambda}(A_2)}\leq C
$
and since $\int_{A_2}\delta(x_1,x',x'') dx_1=0$,  Lemma \ref{lem:trace-sob} implies 
$\Vert \delta\Vert_{H^{1/2}_{\lambda',\lambda''}(\Sigma)}\leq C/\sqrt{\lambda_1}.$
Combined with \eqref{eq:estimethetaH012}, this proves that 
\be\label{eq:estimethetaH12}
\Vert\theta\Vert_{H^{1/2}_{\lambda',\lambda''}(\Sigma)}\leq C/\sqrt{\lambda_1}.
\ee

We are now in position to define the function $\psi$.
Let $\rho\in C^\infty(\mathbb{R_+})$ be such that $\rho(0)=1$ and $\supp(\rho)\subset [0,\frac 12]$.
We define $\psi$ via its partial Fourier coefficients in the variables $(x',x'')$. For $x_1\in]0,1[$ and $\tilde k=(k',k'')\in \Z^{d'+d''}$,  let
$$\hat\psi(x_1,\tilde k)=
\left\{
\begin{array}{c}
\rho\Big(\frac{x_1 }{h\lambda_1}\Big)\hat \theta(\tilde k)\text{ if }\<\tilde k\>\leq h^{-1}\\
\rho\Big(\frac{x_1\langle \tilde k\rangle }{\lambda_1}\Big)\hat \theta(\tilde k)\text{ if }\<\tilde k\>\geq h^{-1}
\end{array}
\right.$$
where for sake of shortness we denote 
$\hat u=\fff_{\lambda',\lambda''} (u)$. Of course, $\psi_{\vert\Sigma}=\theta$ since one has $\hat\psi(0,\tilde k)=\hat\theta(\tilde k)$. 
Moreover, since $\rho$ is supported in $[0,\frac 12]$, 
then $\psi$ is supported in $0< x_1\leq h\lambda_1$. Let us now estimate its 
$L^2$ and $H^1$ norms.
Denoting   $\Vert u\Vert_{L^2(]0,1[\times \Z^{d-1})}^2=\sum_{k\in\Z^{d-1}}\Vert u(.,k)\Vert_{L^2(]0,1[)}^2$,
we have
 \begin{equation*}
 \begin{split}
 \lvert \lvert \psi\rvert \rvert^2_{L^2(A_1)}& = \Vert \hat {\psi}(x_1,\tilde k) \Vert ^2 _{L^2(]0,1[\times\Z^{d-1})}\\
& =\sum_{\<\tilde k\>\leq h^{-1}} \lvert \hat{\theta}(\tilde k)\rvert^2 \int_0^\infty \lvert \rho(\frac{x_1 }{h\lambda_1})\rvert ^2 dx_1+
\sum_{\<\tilde k\>\geq h^{-1}} \lvert \hat{\theta}(\tilde k)\rvert^2 \int_0^\infty\lvert \rho(\frac{x_1\langle \tilde k \rangle }{\lambda_1})\rvert ^2 dx_1\\
& \leq  \Vert \rho\Vert ^2_{L^2}\Big(\sum_{\<\tilde k\>\leq h^{-1}}h\lambda_1\lvert \hat{\theta}(\tilde k)\rvert^2
+\sum_{\<\tilde k\>\geq h^{-1}}\frac{\lambda_1}{\<\tilde k\>}\lvert \hat{\theta}(\tilde k)\rvert^2\Big)\\
&\leq  h \lambda_1 \Vert \rho\Vert ^2_{L^2(\R_+)} \Vert \theta\Vert ^2_{L^2(\Sigma)}\leq Ch^2
\end{split}
\end{equation*}
thanks to \eqref{eq:estimthetaL2}. This proves the second part of \eqref{deux}. To prove the $H^1$ estimate, we  observe that 
\be\label{eq:decompH1norm}
\Vert \psi\Vert ^2_{H^1_{\lambda}(A_1)} = \Vert \psi\Vert ^2_{L^2(A_1)}+\Vert \lambda_1 \partial_1 \hat{\psi}\Vert ^2_{L^2(]0,1[\times\Z^{d-1})}
+\Vert \langle \tilde k\rangle \hat{\psi}\Vert^2_{L^2(]0,1[\times\Z^{d-1})}
\ee
and we estimate separately each term of the right hand side. First, we have
\begin{equation*}
\begin{split}
\Vert \lambda_1 \partial_1 \hat{\psi}\Vert ^2_{L^2(]0,1[\times\Z^{d-1})}
&=\sum_{\<\tilde k\>\leq h^{-1}}h^{-2}
 \lvert \hat{\theta}(\tilde k)\rvert^2 \int_0^\infty \lvert \rho'(\frac{x_1 }{h\lambda_1})\rvert ^2 dx_1\\
 &\phantom{*********}+\sum_{\<\tilde k\>\geq h^{-1}}\<\tilde k\>^2
 \lvert \hat{\theta}(\tilde k)\rvert^2 \int_0^\infty \lvert \rho'(\frac{x_1\langle \tilde k \rangle }{\lambda_1})\rvert ^2 dx_1\\
 &=\Vert \rho'\Vert_{L^2}^2\Big(\frac{\lambda_1}h\sum_{\<\tilde k\>\leq h^{-1}}\lvert \hat{\theta}(\tilde k)\rvert^2+\lambda_1\sum_{\<\tilde k\>\geq h^{-1}}\<\tilde k\>\lvert \hat{\theta}(\tilde k)\rvert^2\Big)\\
 &\leq \Vert \rho'\Vert_{L^2}^2\big(\frac{\lambda_1}h\Vert \theta\Vert_{L^2}^2+\lambda_1\Vert\theta\Vert^2_{H^{\frac 12}_{\lambda',\lambda''}(\Sigma)}\big)\leq C\Vert \rho'\Vert_{L^2}^2
\end{split}
\end{equation*}
thanks to \eqref{eq:estimthetaL2} and \eqref{eq:estimethetaH12}. Let us now estimate the last term in \eqref{eq:decompH1norm}. We have
\begin{equation*}
\begin{split}
\Vert \langle \tilde k\rangle \hat{\psi}\Vert^2_{L^2(]0,1[\times\Z^{d-1})}&=\sum_{\<\tilde k\>\leq h^{-1}}
 \lvert \hat{\theta}(\tilde k)\rvert^2 \<\tilde k\>^2 \int_0^\infty \lvert \rho(\frac{x_1 }{h\lambda_1})\rvert ^2 dx_1\\
 &\phantom{*******}+\sum_{\<\tilde k\>\geq h^{-1}}
 \lvert \hat{\theta}(\tilde k)\rvert^2 \<\tilde k\>^2\int_0^\infty \lvert \rho(\frac{x_1\langle \tilde k \rangle }{\lambda_1})\rvert ^2 dx_1\\
 &=\Vert \rho\Vert_{L^2}^2\Big(h\lambda_1\sum_{\<\tilde k\>\leq h^{-1}}
 \lvert \hat{\theta}(\tilde k)\rvert^2 \<\tilde k\>^2
 +\lambda_1\sum_{\<\tilde k\>\geq h^{-1}}
 \lvert \hat{\theta}(\tilde k)\rvert^2 \<\tilde k\>
 \Big)\\
&\leq \Vert \rho\Vert_{L^2}^2\lambda_1\Vert\theta\Vert^2_{H^{\frac 12}_{\lambda',\lambda''}(\Sigma)}\leq C\Vert \rho\Vert_{L^2}^2
\end{split}
\end{equation*}
thanks again to   \eqref{eq:estimethetaH12}. This achieves to prove that 
$\Vert \psi\Vert ^2_{H^1_{\lambda}(A_1)}=\ooo(1)$. 
\ep

\bibliographystyle{amsplain}
\bibliography{ref_metrocusp}

\end{document}